\documentclass[a4paper]{scrartcl}

\usepackage[a4paper,verbose]{geometry}
\pdfoutput=1

\usepackage[utf8]{inputenc}
\usepackage[T1]{fontenc} %
\usepackage{libertine}

\usepackage{amsfonts}
\usepackage{amssymb}
\usepackage{amsthm}
\usepackage{amsmath}
\usepackage{caption}
\usepackage{etoolbox}
\usepackage{stmaryrd}
\usepackage{ifthen}
\usepackage{suffix}
\usepackage{verbatim} %

\usepackage{xfrac}

\usepackage{bm}
\usepackage{physics}
\usepackage[safe]{tipa} %
\usepackage{mathtools} %
\usepackage{color}
\usepackage{wrapfig}
\usepackage{tikz}
\usepackage{tikzit}
\usetikzlibrary{cd, babel}

\usepackage[framemethod=tikz]{mdframed}

\newcommand*{\secref}[1]{\S\ref{#1}}

\usepackage[backend=biber,style=numeric-comp,sorting=none,date=short,maxcitenames=2]{biblatex}

\usepackage{hyperref}

\newcommand*{\email}[1]{%
  \normalsize\href{mailto:#1}{\texttt{#1}}\par
}

\bibliography{bibliography}

\newcommand*{\citep}[1]{\parencite{#1}}
\newcommand*{\citet}[1]{\textcite{#1}}

\usepackage[capitalize]{cleveref} %

\usepackage{graphicx}
\usepackage{pict2e}

\def\mdash{{\hbox{-}}}

\makeatletter
\newcommand{\adjunction}{\@ifstar\named@adjunction\normal@adjunction}
\newcommand{\normal@adjunction}[4]{%
  #1\colon #2%
  \mathrel{\vcenter{%
    \offinterlineskip\m@th
    \ialign{%
      \hfil$##$\hfil\cr
      \longrightharpoonup\cr
      \noalign{\kern-.3ex}
      \smallbot\cr
      \longleftharpoondown\cr
    }%
  }}%
  #3 \noloc #4%
}
\newcommand{\named@adjunction}[4]{%
  #2%
  \mathrel{\vcenter{%
    \offinterlineskip\m@th
    \ialign{%
      \hfil$##$\hfil\cr
      \scriptstyle#1\cr
      \noalign{\kern.1ex}
      \longrightharpoonup\cr
      \noalign{\kern-.3ex}
      \smallbot\cr
      \longleftharpoondown\cr
      \scriptstyle#4\cr
    }%
  }}%
  #3%
}
\newcommand{\longrightharpoonup}{\relbar\joinrel\rightharpoonup}
\newcommand{\longleftharpoondown}{\leftharpoondown\joinrel\relbar}
\newcommand\noloc{%
  \nobreak
  \mspace{6mu plus 1mu}
  {:}
  \nonscript\mkern-\thinmuskip
  \mathpunct{}
  \mspace{2mu}
}
\newcommand{\smallbot}{%
  \begingroup\setlength\unitlength{.15em}%
  \begin{picture}(1,1)
  \roundcap
  \polyline(0,0)(1,0)
  \polyline(0.5,0)(0.5,1)
  \end{picture}%
  \endgroup
}
\makeatother

\let\op=\relax

\def\op{\ensuremath{^{\,\mathrm{op}}}}

\newcommand{\nn}{{\mathbb{N}}}
\newcommand{\rr}{{\mathbb{R}}}

\newcommand{\Cat}[1]{\mathbf{#1}}
\newcommand{\cat}[1]{\mathcal{#1}}
\newcommand{\Fun}[1]{\mathsf{#1}}

\newcommand{\Kl}{\mathcal{K}\mspace{-2mu}\ell}

\DeclareMathOperator*{\E}{\mathbb{E}}

\renewcommand{\d}{\mathrm{d}}

\newcommand{\para}{\Cat{Para}}

\newcommand{\Da}{{\mathcal{D}}}

\newcommand{\Fa}{{\mathcal{F}}}

\newcommand{\Qa}{{\mathcal{Q}}}

\newcommand{\Giry}{{\mathcal{G}}}

\DeclareMathOperator{\id}{\mathsf{id}}

\DeclareMathOperator{\Set}{\Cat{Set}}

\newcommand{\xto}[1]{\xrightarrow{#1}}

\makeatletter
\newcommand{\mathoverlap}[2]{\mathpalette\mathoverlap@{{#1}{#2}}}
\newcommand{\mathoverlap@}[2]{\mathoverlap@@{#1}#2}
\newcommand{\mathoverlap@@}[3]{\ooalign{$\m@th#1#2$\crcr\hidewidth$\m@th#1#3$\hidewidth}}
\makeatother

\newcommand{\klcirc}{\bullet} %
\newcommand*{\smallklcirc}{\raisebox{0.18ex}{\scalebox{0.66}{$\klcirc$}}}
\newcommand{\klto}{\mathoverlap{\rightarrow}{\smallklcirc\,}}

\newcommand{\xklto}[1]{\mathoverlap{\xrightarrow{#1}}{\smallklcirc\,}}

\def\lenscirc{\baro}
\newcommand{\lensto}{\mathrel{\ooalign{\hfil$\mapstochar\mkern5mu$\hfil\cr$\to$\cr}}}
\newcommand{\xlensto}[1]{\mathoverlap{\xrightarrow{#1}}{\raisebox{0.375ex}{\scalebox{1.0}[0.33]{$|$}}\,}}

\makeatletter
\providecommand*{\xmapstofill@}{%
  \arrowfill@{\mapstochar\relbar}\relbar\rightarrow
}
\providecommand*{\xmapsto}[2][]{%
  \ext@arrow 0395\xmapstofill@{#1}{#2}%
}

\makeatletter
\def\slashedarrowfill@#1#2#3#4#5{%
  $\m@th\thickmuskip0mu\medmuskip\thickmuskip\thinmuskip\thickmuskip
   \relax#5#1\mkern-7mu%
   \cleaders\hbox{$#5\mkern-2mu#2\mkern-2mu$}\hfill
   \mathclap{#3}\mathclap{#2}%
   \cleaders\hbox{$#5\mkern-2mu#2\mkern-2mu$}\hfill
   \mkern-7mu#4$%
}
\def\rightslashedarrowfill@{%
  \slashedarrowfill@\relbar\relbar\mapstochar\rightarrow}
\newcommand\xslashedrightarrow[2][]{%
  \ext@arrow 0055{\rightslashedarrowfill@}{#1}{#2}}
\makeatother

\theoremstyle{definition}
\newtheorem{defn}{Definition}[section]
\newtheorem{notation}[defn]{Notation}
\newtheorem{depict}[defn]{Depiction}
\newtheorem{ex}[defn]{Example}

\newtheorem{rmk}[defn]{Remark}
\newtheorem{obs}[defn]{Observation}

\newtheorem*{rmk*}{Remark}

\newtheorem{prop}[defn]{Proposition}
\newtheorem{prop*}{Proposition}

\newtheorem{lemma}[defn]{Lemma}
\newtheorem{thm}[defn]{Theorem}
\newtheorem{cor}[defn]{Corollary}

\newtheorem*{thm*}{Theorem}
\newtheorem*{cor*}{Corollary}

\definecolor{darkblue}{rgb}{0,0,0.7}

\usepackage{authblk}

\author{Toby St. Clere Smithe}
\affil{University of Oxford \\ \& \\ Topos Institute \\ \email{toby@topos.institute}}

\hypersetup{
  pdftitle={Active Inference I},
  pdfauthor={Toby St Clere Smithe}
}

\tikzstyle{xshiftu}=[shift = {(#1, 0)}]
\tikzstyle{yshiftu}=[shift = {(0, #1)}]

\tikzstyle{dot}=[inner sep=0.25mm,minimum width=1mm,minimum height=1mm,draw,shape=circle,text depth=-0.2mm]
\tikzstyle{white dot}=[dot,fill=white, draw=black]
\tikzstyle{action}=[dot,fill=white,scale=0.667,inner sep=0.5mm]
\tikzstyle{copier}=[dot,fill=white,scale=2.0]
\tikzstyle{black copier}=[dot,fill=black,scale=2.0]

\tikzstyle{box}=[fill=white, draw=black, shape=rectangle]
\tikzstyle{medium box}=[fill=white, draw=black, shape=rectangle, minimum width=1.5cm, minimum height=0.66cm]
\tikzstyle{arrow box}=[fill=white, draw, shape=rectangle,minimum height=5mm,yshift=-0.5mm,minimum width=5mm]

\tikzstyle{effect}=[regular polygon, regular polygon sides=3,draw]
\tikzstyle{state0}=[regular polygon, regular polygon sides=3,draw,shape border rotate=0]
\tikzstyle{state90}=[regular polygon, regular polygon sides=3,draw,shape border rotate=90]
\tikzstyle{state180}=[regular polygon, regular polygon sides=3,draw,shape border rotate=180]
\tikzstyle{state270}=[regular polygon, regular polygon sides=3,draw,shape border rotate=270]

\tikzstyle{scalar}=[diamond,draw,inner sep=1pt]

\tikzstyle{discarder}=[my ground,draw,inner sep=0pt,minimum width=4.2pt,minimum height=11.2pt,anchor=input,rotate=90]
\tikzstyle{discarder0}=[my ground,draw,inner sep=0pt,minimum width=4.2pt,minimum height=11.2pt,anchor=input,rotate=0]

\tikzstyle{pointy1}=[->]
\tikzstyle{midpoint1}=[-, {postaction={decorate,decoration={markings, mark=at position .5 with {\arrow{>}}}}}]
\tikzstyle{midpointy1pointy}=[->, {postaction={decorate,decoration={markings, mark=at position .5 with {\arrow{>}}}}}]
\tikzstyle{dashed1}=[-, dashed]
\tikzstyle{dotted1}=[-, dotted]
\tikzstyle{dash-pointy}=[->, dashed]

\input{strings.tikzdefs}

\ifexternalizetikz\tikzexternaldisable\fi
\newsavebox\sbground
\savebox\sbground{%
  \begin{tikzpicture}[baseline=0pt]
    \draw (0,-.1ex) to (0,.85ex)
    node[ground IEC,draw,anchor=input,inner sep=0pt,
    minimum width=3.15pt,minimum height=8.4pt,rotate=90] {};
  \end{tikzpicture}%
}
\newcommand{\ground}{\mathord{\usebox\sbground}}

\newsavebox\sbcopier
\savebox\sbcopier{%
  \begin{tikzpicture}[baseline=0pt]
    \node[copier,scale=0.7] (a) at (0,3.8pt) {};
    \draw (a) -- +(-90:.21);
    \draw (a) -- +(45:.21);
    \draw (a) -- +(135:.21);
  \end{tikzpicture}}

\ifexternalizetikz\tikzexternalenable\fi

\newsavebox\bsbcopier
\savebox\bsbcopier{%
  \begin{tikzpicture}[baseline=0pt]
    \node[black copier,scale=0.7] (a) at (0,3.8pt) {};
    \draw (a) -- +(-90:.21);
    \draw (a) -- +(45:.21);
    \draw (a) -- +(135:.21);
  \end{tikzpicture}}
\newcommand{\bcopier}{\mathord{\usebox\bsbcopier}}
\ifexternalizetikz\tikzexternalenable\fi

\usepackage{mathabx}

\newcommand{\longhookrightarrow}{\lhook\joinrel\longrightarrow}

\usepackage{quiver}

\newcommand{\BLens}[1]{\Cat{BayesLens}_{\cat{#1}}} %

\newcommand{\BLCtx}[1]{\overline{\BLens{#1}}}    %
\newcommand{\BLCCtx}[1]{\widetilde{\BLens{#1}}}
\def\ctx{\Fun{Ctx}}
\def\cctx{\mathbb{C}\Fun{tx}}

\newcommand{\efb}[3]{\llparenthesis\,{#1}\,|\,{#2}\,|\,{#3}\,\rrparenthesis}

\newcommand{\SGame}[2][{}]{{}_{#1}\Cat{SGame}_{#2}} %
\newcommand{\SimpSGame}[1]{\Cat{SimpSGame}_{#1}}   %

\newcommand{\vlens}[2]{\ensuremath{\begin{pmatrix} #1 \\ #2 \end{pmatrix}}}
\WithSuffix\newcommand\vlens*[2]{\ensuremath{\Bigl(\negthinspace\begin{smallmatrix}#1\\#2\end{smallmatrix}\Bigr)}}

\date{\today}
\title{Compositional Active Inference I: \\ Bayesian Lenses. Statistical Games.}

\begin{document}

\maketitle

\begin{abstract}
  We introduce the concepts of \textit{Bayesian lens}, characterizing the
  bidirectional structure of exact Bayesian inference, and \textit{statistical
    game}, formalizing the optimization objectives of approximate inference
  problems. We prove that Bayesian inversions compose according to the
  compositional lens pattern, and exemplify statistical games with a number of
  classic statistical concepts, from maximum likelihood estimation to
  generalized variational Bayesian methods. This paper is the first in a series
  laying the foundations for a compositional account of the theory of active
  inference, and we therefore pay particular attention to statistical games with
  a free-energy objective.
\end{abstract}

\section{Introduction}

Those systems that we might classify as `cybernetic', `adaptive', or `alive' all
display a fundamental property: they resist perturbations that would push them
away from their goals or render their existence unsustainable. In order to do
so, such systems are somehow able to sense their current state of affairs
(through \textit{perception}) and respond appropriately (through
\textit{action}). In the series of papers of which this is the first part, we
seek to supply new compositional foundations for a theory of \textit{active
  inference} adequate to describe such systems, with a particular focus on the
framework that has come to be known in the compositional neuroscience and
artificial life communities as the \textit{free energy principle}
\parencite{Buckley2017free}, whose structures we seek to make precise.

A central feature of active inference is the use of the statistical procedure
called \textit{Bayesian inference}, which supplies a recipe by which a system
might invert a statistical model (say, of how causes generate observations) in
order to form beliefs about the causes of observed data. It is easy to see how
such a process of inferring causes could be understood as a process of
perception, but the central dogma of active inference is that both perception
and action can be rendered as problems of Bayesian inference, with action being
`dual' to perception: instead of changing its \textit{internal} state (its
beliefs about causes) to match its observations better, a system might act to
change the \textit{external} state (the state of the world) so that the
observations that it expects or desires obtain. In the free energy framework,
both perception and action emerge through the optimization of a single quantity,
the \textit{free energy}.

Such processes of optimization, and perception and action more generally, are
inherently dynamical processes. In this first paper of the series, we put the
dynamics temporarily aside, and lay the statistical foundations, characterizing
the compositional structure of the \textit{generative models} instantiated by
cybernetic systems and the algebra of their inversion. This algebra is
formalized by our concept of \textit{Bayesian lens}, which we introduce to
characterize the inherently bidirectional structure of Bayesian inversion,
drawing on the `lens' pattern that structures bidirectional systems from
economic games \parencite{Ghani2016Compositional}, to databases
\parencite{Bohannon2006Relational}, and machine-learners
\parencite{Fong2019Lenses}.

This lens structure serves more than an organizing purpose: we prove that the
inversion of a composite or `hierarchical' statistical model is equivalently
given (up to almost-equality) by the lens composition rule. This means that
cybernetic systems embodying complex composite models can simply invert each
component factor of their models and then combine these inversions, in order to
obtain an inversion of the whole. In turn, this explains the observation that
hierarchical systems in the brain (such as much of the visual cortex) can be
explained as a composite of `local' circuits each performing a form of
approximate Bayesian inference call predictive coding
\parencite{Bastos2012Canonical}.

Having established the structures required to state and prove that ``Bayesian
updates compose optically''\footnote{Although we do not use the heavier
machinery of `optics' here; for that, see our preprint
\parencite{Smithe2020Bayesian} and Remark \ref{rmk:optics} below.}, we formalize
the ``algebra of statistical inference problems'' as categories of
\textit{statistical games}. These `games' consist of a lens paired with a
contextual \textit{fitness function}, which define the quantities that we often
think of cybernetic systems as optimizing, and where the `context' formalizes
the system's interaction with its environment. In this development, we draw much
inspiration from compositional game theory
\parencite{Ghani2016Compositional,Bolt2019Bayesian}. We exemplify these
statistical games with a range of examples from maximum likelihood estimation to
generalized variational Bayesian methods.

This paper is the first of a series of papers. The next instalment introduces
the structures necessary to supply statistical games with ``dynamical
semantics'', and thus breathe some life into those systems that perform
approximate inference. A subsequent paper will then explain how such systems can
perform action, and thereby affect the worlds that they inhabit.

\paragraph{Overview of this paper}

We begin in \secref{sec:bg-comp-prob} by introducing the basics of compositional
probability theory. In \secref{sec:bg-lens}, we introduce the structures
necessary to describe the lens pattern, and recall that in `nice' situations,
the resulting categories of lenses are monoidal. In \secref{sec:bg-para}, we
introduce the \(\para\) construction, which has been proposed
\parencite{Capucci2021Towards} as foundational for categorical cybernetics, and
often plays an important role for us, too. Then, in \secref{sec:buco}, we define
\textit{Bayesian lenses} and prove the theorem that Bayesian inversions compose
according to the lens pattern; we also give a detailed description of the
low-dimensional structure of parameterized Bayesian lenses
(\secref{sec:buco-para-blens}). In \secref{sec:sgame}, we define contexts for
Bayesian lenses, fitness functions, and the resulting monoidal categories of
statistical games. Finally, in \secref{sec:examples}, we give a number of
examples.

\paragraph{Acknowledgements}

This series of papers is the result of a number of interactions within and
around the applied category theory and active inference communities. We thank
the reviewers and organizers of the ACT conferences in 2020, 2021, and 2022, as
well as those of the SYCO series of symposia. We thank the Topos Institute and
Foundational Questions Institute for support and funding. And we particularly
thank the following individuals (in no particular order) for stimulating
discussions, comments, and encouragement: Samson Abramsky; Matteo Capucci; Bob
Coecke; Lance da Costa; Brendan Fong; Karl Friston; Bruno Gavranović; Neil
Ghani; Jules Hedges; Johannes Kleiner; Tim Hosgood; Sophie Libkind; David Jaz
Myers; Valeria de Paiva; Evan Patterson; Maxwell Ramstead; Dalton Sakthivadivel;
Brandon Shapiro; David Spivak; Sam Staton; Sean Tull; Vincent Wang.

\section{Mathematical background for statistical games}

\subsection{Compositional probability, concretely and abstractly} \label{sec:bg-comp-prob}

In order to define Bayesian lenses and statistical games and prove some basic
results about them, we will work at a high level of abstraction; then, to
exemplify them with applications, we will need to work more concretely.  In both
instances, our basic categorical setting will be a \textit{copy-delete} or
\textit{Markov category}, whose morphisms we will call \textit{stochastic
  channels} or \textit{Markov kernels} and which behave like functions with
uncertain outputs; for us, these model the processes by which observational data
are (believed to be) generated by processes in the world, and composite channels
model composite (sequences of) processes.

We will typically be interested in applications where the sample spaces are
continuous and where the probability measures may have infinite support. A
stochastic channel \(c : X \klto Y\) will be something like a function taking
values in probability measures over a space, but in applications one will often
fix a reference measure and then work with a density function \(p_c : X \times Y
\to [0, 1]\) representing the channel \(c\) with respect to that measure. We
begin this section by introducing the category \(\Cat{sfKrn}\) of
\textit{s-finite kernels}, where these concepts have precise meanings, before
generalizing graphically to the abstract setting of copy-delete categories.

\subsubsection{S-finite kernels}

We sketch the basic structure of the category \(\Cat{sfKrn}\) of s-finite
kernels between measurable spaces, and refer the reader to
\citet{Cho2017Disintegration,Staton2017Commutative} for elaboration of the
details.

\begin{defn}
  Suppose \((X, \Sigma_X)\) and \((Y, \Sigma_Y)\) are measurable spaces, with
  \(X,Y\) sets and \(\Sigma_X,\Sigma_Y\) the corresponding \(\sigma\)-algebras.
  A \emph{kernel} \(k\) from \(X\) to \(Y\) is a function \(k : X \times
  \Sigma_Y \to [0, \infty]\) satisfying the conditions:
  \begin{itemize}
  \item for all \(x \in X\), \(k(x, -) : \Sigma_Y \to [0, \infty]\) is a measure; and
  \item for all \(B \in \Sigma_Y\), \(k(-, B) : X \to [0, \infty]\) is measurable.
  \end{itemize}
  A kernel \(k : X \times \Sigma_Y \to [0, \infty]\) is \emph{finite} if there
  exists some \(r \in [0, \infty)\) such that, for all \(x \in X\), \(k(x, Y)
  \leq r\). And \(k\) is \emph{s-finite} if it is the sum of at most countably
  many finite kernels \(k_n\), \(k = \sum_{n : \nn} k_n\).
\end{defn}

\begin{prop}
  Measurable spaces and s-finite kernels \((X, \Sigma_X) \klto (Y, \Sigma_Y)\)
  between them form a category, denoted \(\Cat{sfKrn}\); note that often we will
  just write \(X\) for a measurable space, leaving the \(\sigma\text{-algebra }
  \Sigma_X\) implicit. Identity morphisms \(\id_X : X \klto X\) in
  \(\Cat{sfKrn}\) are Dirac kernels \(\delta_X : X \times \Sigma_X \to [0,
    \infty] := x \times A \mapsto 1\) iff \(x \in A\) and 0
  otherwise. Composition is given by a Chapman-Kolmogorov equation: suppose \(c
  : X \klto Y\) and \(d : Y \klto Z\). Then
  \[
  d \klcirc c : X \times \Sigma_Z \to [0, \infty]
  := x \times C \mapsto \int_{y:Y} d(C|y) \, c(\d y| x)
  \]
  where we have used `conditional probability' notation \(d(C|y) := d(y, C)\).
\end{prop}

\begin{rmk}
  In the scientific literature, one often encounters the term
  \textit{conditional probability distribution}, which can almost always be
  interpreted as indicating a probability kernel: we can think of a probability
  kernel \(c : X \klto Y\) as a function that emits a probability distribution
  \(c(x) : \Sigma_Y \to [0, \infty]\) over Y for each \(x : X\). We can then see
  \(c(x)\) as a probability distribution \textit{conditional on} the choice or
  observation of \(x : X\), hence the notation \(c(-|x)\). One reads this
  notation \(-|x\) as ``\(-\) \textit{given} \(x\)''.
\end{rmk}

\begin{prop}
  There is a monoidal structure \((\otimes, 1)\) on \(\Cat{sfKrn}\), the unit of
  which is the singleton set \(1\) with its trivial sigma-algebra. On objects,
  \(X \otimes Y\) is the Cartesian product \((X \times Y, \Sigma_{X \times Y})\)
  of measurable spaces, with the product sigma-algebra. On morphisms, \(f
  \otimes g : X \otimes Y \klto A \otimes B\) is given by
  \[
  f \otimes g : (X \times Y) \times \Sigma_{A \times B}
  := (x \times y) \times E \mapsto \int_{a:A} \int_{b:B} \delta_{A \otimes B}(E|x, y) \, f(\d a|x) \, g(\d b|y)
  \]
  where, as above, \(\delta_{A \otimes B}(E|a, b) = 1\) iff \((a, b) \in E\) and
  0 otherwise. Note that \((f \otimes g)(E|x, y) = (g \otimes f)(E|y, x)\) for
  all s-finite kernels (and all \(E\), \(x\) and \(y\)), by the Fubini-Tonelli
  theorem for s-finite measures
  \citep{Cho2017Disintegration,Staton2017Commutative}, and so \(\otimes\) is
  symmetric on \(\Cat{sfKrn}\).
\end{prop}

\begin{rmk}[States and effects]
  We will call kernels with domain \(1\) \textit{states}. A kernel \(1 \klto X\)
  is equivalently a function \(\Sigma_X \to [0, \infty]\), which is simply a
  (possibly improper) measure on the space \((X, \Sigma_X)\). Occasionally we
  will say \textit{distribution} to mean `state'. We call a state on a product
  space, such as \(\omega : 1 \klto X \otimes Y\), a \textit{joint state} or
  \textit{joint distribution}.

  Dually, kernels with codomain \(1\) will be called \textit{effects}. Note that
  although \(1\) is the unit of the preceding monoidal structure, this unit is
  not terminal in \(\Cat{sfKrn}\): s-finite kernels \(X \klto 1\) are
  equivalently measurable functions \(X \to [0, \infty]\), and there are of
  course many nontrivial examples; of which density functions will form an
  important class.
\end{rmk}

\begin{defn}[Probability kernel, probability measure, probability space]
  If a kernel \(k : X \times \Sigma_Y \to [0,\infty]\) satisfies the additional
  conditions that it takes values in the unit interval \([0,1]\) and, for all
  \(x : X\), \(k(Y) = 1\), then we call \(k\) a \textit{probability kernel}. If
  \(k\) is a state (\textit{i.e.}, \(X = 1\)), then we call it a
  \textit{probability measure}. A \textit{probability space} is a pair
  \((\Omega, \pi)\) of a measurable space \(\Omega\) with a probabilty measure
  \(\pi : 1 \klto \Omega\).
\end{defn}

\begin{rmk}[Giry monad]
  Probability measures \(1 \klto X\) on each \(X\) form the points of a space
  which we will denote \(\Giry X\). This space can be equipped with a canonical
  \(\sigma\)-algebra, making \(\Giry\) into a functor \(\Cat{Meas} \to
  \Cat{Meas}\) which acts on each measurable function \(f : X \to Y\) by
  returning its \textit{pushforward} \(f_* : \Giry X \to \Giry Y\), defined by
  \(f_* \nu : \Sigma_Y \to [0, 1] : B \mapsto \nu(f^{-1}(B))\). \(\Giry\) can in
  turn be equipped with a monad structure \((\mu, \eta)\). Every measurable
  function \(X \to \Giry Y\) corresponds to a probability kernel \(X \klto Y\),
  making the Kleisli category \(\Kl(\Giry)\) a subcategory of \(\Cat{sfKrn}\).
  Composition of probability kernels \(X \xklto{f} Y \xklto{g} Z\) corresponds
  accordingly to Kleisli composition \(X \xto{f} \Giry Y \xto{\Giry g} \Giry
  \Giry Z \xto {\mu_Z} \Giry Z\), and the components of the monad unit
  correspond to the Dirac (identity) kernels \(\delta_X : X \klto X\).  We refer
  the reader to \textcite[{\S 4}]{Fritz2019synthetic} and the references therein
  for elaboration.  The slice \(1/\Kl(\Giry)\) of \(\Kl(\Giry)\) under \(1\) is
  the category of probability spaces and measure-preserving kernels between
  them, called \(\Cat{ProbStoch}\) by \textcite{Fritz2019synthetic}.
\end{rmk}

\begin{rmk}[Convex spaces and expectations] \label{rmk:convex}
  A \textit{convex space} is an algebra of the Giry monad; that is, a space
  \(X\) equipped with a measurable function \(\Giry X \to X\) called the
  algebra evaluation or \textit{expected value}. Each measurable function \(f :
  X \to X\) induces an expected value \(\E_f : \Giry X \to X\) defined as
  \[
  \mathbb{E}_f(\pi) := \int_{x:X} f(x) \, \pi(\d x) \, .
  \]
  We will typically instead write
  \[
  \E_{x \sim \pi} \, [f] := \mathbb{E}_f(\pi)
  \]
  where the notation \(x \sim \pi\) should be read as ``\(x\) distributed
  according to \(\pi\)''. More generally, we have an operator \(\E :
  \Cat{Meas}(\Omega, X) \times \Giry \Omega \to X\) defined similarly by
  \[
  \E_{\omega \sim \pi} \, [p] := \int_{\omega:\Omega} p(\omega) \, \pi(\d \omega)
  \]
  where \(p : \Omega \to X\) and \(\pi : \Giry \Omega\).  Note that this
  subsumes the case where \(p\) is an \(X\)-valued random variable defined on a
  probability space \((\Omega, \pi)\). Commonly, we will have \(X = \rr\) or \(X
  = [0, \infty]\).
\end{rmk}

\begin{rmk}[Effects and validities]
  A special case of the preceding expectation operator occurs when \(X = [0,
    \infty]\). In this case, maps \(\Omega \to [0,\infty]\) are of course
  effects \(\Omega \klto 1\) in \(\Cat{sfKrn}\), and for each state \(\pi : 1
  \klto \Omega\) and effect \(p : \Omega \klto 1\), the expectation operator
  simply computes the composte \(p \klcirc \omega\).  That is, we have in this
  case \(\E_{\omega\sim\pi}[p] = p \klcirc \omega\), and we might then call
  \(p\) a \textit{predicate} on \(\Omega\) and the expectation
  \(\E_{\omega\sim\pi}[p]\) the \textit{validity} of \(p\) in the state
  \(\pi\). We refer the reader to \textcite[\S 5]{Cho2015Introduction} for more
  on this perspective, where the validity is written \(\pi \models p\).
\end{rmk}

\begin{obs} \label{obs:ground}
  Each space \(X\) in \(\Cat{sfKrn}\) is equipped with a canonical effect
  \(\ground_X : X \klto 1\), the constant effect \(x \mapsto 1\). We denote this
  family of effects by the `ground' symbol \(\ground\) and call the components
  \textit{discarding maps} with the intuition that they act by `discarding'
  information (``wiring to ground''). Another way to characterize probability
  kernels, and thus \(\Kl(\Giry)\), is that they make discarding natural,
  satisfying \(\ground \klcirc c = \ground\). We call such channels or processes
  \textit{causal}: they cannot affect the outcomes of processes earlier in a
  sequence of composites.
\end{obs}

Not only can we discard information in \(\Cat{sfKrn}\), but we can also
\textit{copy} it. In other words, each object is equipped with a canonical
copy-delete structure, making it into a comonoid, with which we can duplicate
and discard states of the corresponding type. In the terms of
\textcite{Fong2019Supplying}, \(\Cat{sfKrn}\) supplies comonoids.

\begin{prop}[\(\Cat{sfKrn}\) supplies comonoids]
  Each object \(X\) is equipped with a canonical comonoid structure
  \((\bcopier_X, \ground_X)\), with \(\bcopier_X : X \klto X \otimes X\) and
  \(\ground_X : X \klto 1\) satisfying the usual comonoid laws. Discarding is
  given by the family of effects \(\ground_X : X \to [0, \infty] := x \mapsto
  1\), and copying is Dirac-like: \(\bcopier_X : X \times \Sigma_{X \times X} :=
  x \times E \mapsto 1\) iff \((x, x) \in E\) and 0 otherwise.
\end{prop}

Discarding part of a joint state gives us \textit{marginals} (or ``marginal
distributions'').

\begin{defn} \label{def:marginals}
  Given a joint distribution \(\omega : 1 \klto X \otimes Y\), we call
  \(\omega_X := (\id_X \otimes \ground_Y) \klcirc \omega : 1 \klto X\) and
  \(\omega_Y := (\ground_X \otimes \id_Y) \klcirc \omega : 1 \klto Y\) the
  \textit{marginals} of \(\omega\). We define projection (or
  \textit{marginalization}) operators by \(\mathsf{proj}_X := \id_X \otimes
  \ground_Y : X\otimes Y \to X\otimes 1 \xto{\sim} X\) and \(\mathsf{proj}_Y :=
  \ground_X \otimes \id_Y : X\otimes Y \to 1\otimes Y \xto{\sim} Y\).
\end{defn}

In this work, we are interested in the problem of inverting stochastic channels
\(X \klto Y\) in order to obtain channels \(Y \klto X\), and we are particularly
interested in what is known as \textit{Bayesian} inversion. As we will see, the
Bayesian inversion of a channel \(c : X \klto Y\) is determined in conjunction
with a state \(\pi : 1 \klto X\).

\begin{defn} \label{def:gen-model}
  We call the pairing \((\pi, c)\) of a state \(\pi : 1 \klto X\) with a channel
  \(c : X \klto Y\) a \textit{generative model} \(X \klto Y\). It induces a
  joint distribution \(\omega_{(\pi, c)} := (\id_X\otimes\, c) \klcirc
  \bcopier_X \klcirc \pi : 1 \klto X \otimes Y\). The marginals of
  \(\omega_{(\pi, c)}\) are \(\pi\) and \(c\klcirc\pi\).
\end{defn}

In the informal scientific literature, the Bayesian inversion of a channel \(X
\klto Y\) (typically called the `likelihood') with respect to a state \(1 \klto
X\) (typically called the `prior') is often written as the expression
\begin{equation} \label{eq:bayes-density-informal}
  p(x|y)
  = \frac{p(y|x) \, p(x)}{p(y)}
  = \frac{p(y|x) \, p(x)}{\int_{x':X} p(y|x') \, p(x') \, \d x'} \, ,
\end{equation}
but this expression is very ill-defined: what is \(p(y|x)\), and how does it
relate to a channel \(c : X \klto Y\)? Why are the clearly different terms
\(p(x|y), p(y|x), p(x), p(y)\) all written with the same symbol \(p\)?

To answer these questions and clarify such expressions, we use density
functions.

\begin{defn}[Density functions] \label{def:density-func}
  We will say that a kernel \(c : X \klto Y\) is \textit{represented by the
    effect} \(p_c : X \otimes Y \klto 1\) \textit{with respect to the state}
  \(\mu : 1 \klto Y\) when
  \[
  c : X \times \Sigma_Y \to [0, \infty]
  := x \times B \mapsto \int_{y:B} \mu(\d y) \, p_c(y | x) .
  \]
  We call the corresponding function \(p_c : X \times Y \to [0, \infty]\) a
  \textit{density function} for \(c\).  Note that we also use conditional
  probability notation for density functions, and so \(p_c(y|x) := p_c(x, y)\).
\end{defn}

\begin{rmk}
  When a channel \(c\) is associated with a density function, we will adopt the
  convention of naming the density function \(p_c\); that is, with a subscript
  indicating the corresponding channel.  In this way, we can rewrite Equation
  \eqref{eq:bayes-density-informal} as
  \begin{equation} \label{eq:bayes-density-informal-2}
  p_{c^\dag_\pi}(x|y)
  = \frac{p_c(y|x) \, p_\pi(x)}{p_{c\klcirc\pi}(y)}
  = \frac{p_c(y|x) \, p_\pi(x)}{\int_{x':X} p_c(y|x') \, p_\pi(x') \, \d x'} \, .
  \end{equation}
\end{rmk}

\begin{rmk}
  We will also adopt the convention of denoting a Bayesian inversion of the
  channel \(c\) with respect to the state \(\pi\) by the symbol
  \(c^\dag_\pi\). We adopt this symbol because it is known that Bayesian
  inversion induces a `dagger' functor \parencite{Karvonen2019Way} on (a
  quotient of) the category \(1/\Kl(\Giry)\) of probability spaces and
  measure-preserving functions \parencite[Remark 13.9]{Fritz2019synthetic}; when
  we `forget' the measures associated with the probability spaces---which form
  the `priors' for the inversions---then we have to explicitly incorporate them
  into the structure, which we indicate with the subscript \(\pi\) in
  \(c^\dag_\pi\).
\end{rmk}

We wrote ``\textit{a} Bayesian inversion'' in the preceding remark since in a
general measurable setting Bayesian inversions need not always exist
\citep{Stoyanov2014Counterexamples}, and when they do they may only be unique up
to almost-equality.

\begin{defn}[Almost-equality]
  Given a state \(\pi : 1 \klto X\), we say that two parallel channels \(c,d : X
  \klto Y\) are \(\pi\)\textit{-almost-equal}, denoted \(c \overset{\pi}{\sim}
  d\), if the joint distributions of the two generative models \((\pi, c)\) and
  \((\pi, d)\) are equal; that is, if \((\id\otimes\, c) \klcirc \bcopier
  \klcirc \pi = (\id\otimes\, d) \klcirc \bcopier \klcirc \pi\).
\end{defn}

\begin{rmk}
  Like the notion of joint distribution for a generative model, the meaning of
  the definition of \(\pi\)-almost equality---that the induced joint states are
  equal---will be rendered clearer in the graphical calculus.
\end{rmk}

We are now in a position to define Bayesian inversions for channels in
\(\Cat{sfKrn}\), although we leave the abstract definition satisfied by the
following until Definition \ref{def:admit-bayes}.

\begin{prop}[{\textcite[Example 8.4]{Cho2017Disintegration}}]
  Suppose \((\pi, c)\) is a generative model \(X \klto Y\) in \(\Cat{sfKrn}\),
  where \(c\) is represented by the effect \(p_c\) with respect to the state
  \(\mu : 1 \klto X\). Then, when it exists, the channel \(c^\dag_\pi : Y \klto
  X\) defined as follows is a Bayesian inversion of \(c\) with respect to
  \(\pi\):
  \begin{equation*} %
    \begin{aligned}
      c^\dag_\pi : Y \times \Sigma_X \to [0, \infty]
      := y \times A & \mapsto \left( \int_{x:A} \pi(\d x) \, p_c(y|x) \right) p_c^{-1}(y) \\
      &= p_c^{-1}(y) \int_{x:A} p_c(y|x) \, \pi(\d x) ,
    \end{aligned}
  \end{equation*}
  where \(p_c^{-1} : Y \klto I\) is given up to \(\mu\text{-almost-equality}\) by
  \[
  p_c^{-1} : Y \to [0, \infty]
  := y \mapsto \left( \int_{x:X} p_c(y|x) \, \pi(\d x) \right)^{-1} \, .
  \]
\end{prop}

Note that from the preceding proposition we recover the informal form of Bayes'
rule (Equation \eqref{eq:bayes-density-informal-2}). Suppose \(\pi\) is itself
represented by a density function \(p_\pi\) with respect to the Lebesgue measure
\(\d x\). Then
\[
c^\dag_\pi (A|y) = \int_{x:A} \, \frac{p_c(y|x) \, p_\pi(x)}{\int_{x':X} \, p_c(y|x') \, p_\pi(x') \, \d x'} \; \d x.
\]

\subsubsection{Copy-delete categories and their graphical calculus} \label{sec:cd-cat}

While most of our examples and applications will found in \(\Cat{sfKrn}\), most
of our definitions and results hold more generally, and it is in such more
general terms that they are most naturally expressed. Our main language will be
that of categories like \(\Cat{sfKrn}\) in which information can be transformed,
copied, and deleted.

\begin{defn}[{\textcite[Def. 2.2]{Cho2017Disintegration}}] \label{def:cd-cat}
  A \textit{copy-delete category} is a symmetric monoidal category \((\cat{C},
  \otimes, I)\) in which every object \(X\) is supplied with a commutative
  comonoid structure \((\bcopier_X, \ground_X)\) compatible with the monoidal
  structure of \((\otimes, I)\). An \textit{affine} copy-delete category, or
  \textit{Markov category} \citep{Fritz2019synthetic}, is a copy-delete category
  in which every channel \(c\) is causal in the sense that \(\ground \klcirc c =
  \ground\). Equivalently, a Markov category is a copy-delete category in which
  the monoidal unit \(I\) is the terminal object.
\end{defn}

\begin{ex}
  \(\Cat{sfKrn}\) is a copy-delete category, while \(\Kl(\Giry)\) is a Markov
  category.
\end{ex}

Monoidal categories, and (co)monoids within them, admit a formal graphical
calculus that substantially simplifies many calculations involving complex
morphisms: proofs of many equalities reduce to visual demonstrations of isotopy,
and structural morphisms such as the symmetry of the monoidal product acquire
intuitive topological depictions. We make substantial use of this calculus
below, and summarize its features here. For more details, see
\textcite[{\S}2]{Cho2017Disintegration} or \textcite[{\S}2]{Fritz2019synthetic}
or the references cited therein.

\begin{depict}[Basic conventions]
  String diagrams in this paper will be read vertically, with information
  flowing upwards (from bottom to top); in later parts, we will have diagrams
  oriented left-to-right. Sequential composition is represented by connecting
  strings together; and parallel composition \(\otimes\) by placing diagrams
  adjacent to one another.  This way, \(c : X \klto Y\), \(\id_X : X \klto X\),
  \(d \klcirc c : X \xklto{c} Y \xklto{d} Z\), and \(f \otimes g : X \otimes Y
  \klto A \otimes B\) are depicted respectively as:
  
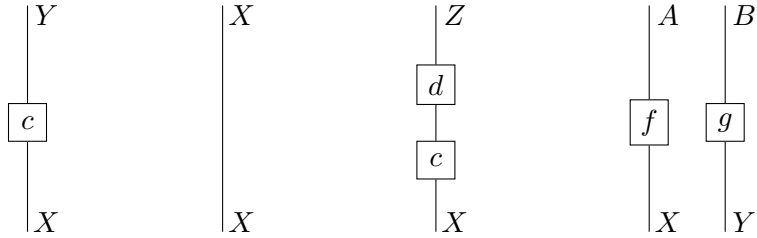
\[
  \hspace{0.125\linewidth} \begin{tikzpicture}
	\begin{pgfonlayer}{nodelayer}
		\node [style=none] (0) at (0, -3) {};
		\node [style=arrow box] (1) at (0, 0) {$c$};
		\node [style=none] (2) at (0, 3) {};
		\node [style=none] (3) at (0.5, -2.75) {$X$};
		\node [style=none] (4) at (0.5, 2.75) {$Y$};
	\end{pgfonlayer}
	\begin{pgfonlayer}{edgelayer}
		\draw (0.center) to (1);
		\draw (1) to (2.center);
	\end{pgfonlayer}
\end{tikzpicture}
  \hspace{0.125\linewidth} \begin{tikzpicture}
	\begin{pgfonlayer}{nodelayer}
		\node [style=none] (0) at (0, -3) {};
		\node [style=none] (2) at (0, 3) {};
		\node [style=none] (3) at (0.5, -2.75) {$X$};
		\node [style=none] (4) at (0.5, 2.75) {$X$};
	\end{pgfonlayer}
	\begin{pgfonlayer}{edgelayer}
		\draw (0.center) to (2.center);
	\end{pgfonlayer}
\end{tikzpicture}

  \hspace{0.125\linewidth} \begin{tikzpicture}
	\begin{pgfonlayer}{nodelayer}
		\node [style=arrow box] (1) at (0, 1) {$d$};
		\node [style=none] (2) at (0, 3) {};
		\node [style=none] (4) at (0.5, 2.75) {$Z$};
		\node [style=none] (5) at (0, -3) {};
		\node [style=arrow box] (6) at (0, -1) {$c$};
		\node [style=none] (8) at (0.5, -2.75) {$X$};
	\end{pgfonlayer}
	\begin{pgfonlayer}{edgelayer}
		\draw (1) to (2.center);
		\draw (5.center) to (6);
		\draw (6) to (1);
	\end{pgfonlayer}
\end{tikzpicture}
  \hspace{0.125\linewidth} \begin{tikzpicture}
	\begin{pgfonlayer}{nodelayer}
		\node [style=none] (0) at (-1, -3) {};
		\node [style=arrow box] (1) at (-1, 0) {$f$};
		\node [style=none] (2) at (-1, 3) {};
		\node [style=none] (3) at (-0.5, -2.75) {$X$};
		\node [style=none] (4) at (-0.5, 2.75) {$A$};
		\node [style=none] (5) at (1, -3) {};
		\node [style=arrow box] (6) at (1, 0) {$g$};
		\node [style=none] (7) at (1, 3) {};
		\node [style=none] (8) at (1.5, -2.75) {$Y$};
		\node [style=none] (9) at (1.5, 2.75) {$B$};
	\end{pgfonlayer}
	\begin{pgfonlayer}{edgelayer}
		\draw (0.center) to (1);
		\draw (1) to (2.center);
		\draw (5.center) to (6);
		\draw (6) to (7.center);
	\end{pgfonlayer}
\end{tikzpicture}
  \hspace{0.125\linewidth}
  \]
  We represent (the identity morphism on) the monoidal unit \(I\) as an empty
  diagram: that is, we leave it implicit in the graphical representation.
\end{depict}

\begin{depict}[States and effects]
  States \(\sigma : I \klto X\) and effects \(\eta : X \klto I\) will be
  depicted as follows:
  \begin{gather*}
    \scalebox{0.85}
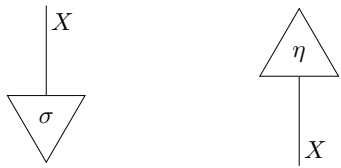
{\begin{tikzpicture}
	\begin{pgfonlayer}{nodelayer}
		\node [style=state180] (0) at (0, -2) {$\sigma$};
		\node [style=none] (1) at (0, 1.5) {};
		\node [style=none] (2) at (0.5, 1) {$X$};
	\end{pgfonlayer}
	\begin{pgfonlayer}{edgelayer}
		\draw (0) to (1.center);
	\end{pgfonlayer}
\end{tikzpicture}
}
    \hspace{0.125\linewidth}
    \scalebox{0.85}{\begin{tikzpicture}
	\begin{pgfonlayer}{nodelayer}
		\node [style=state0] (0) at (0, 0) {$\eta$};
		\node [style=none] (1) at (0, -3.5) {};
		\node [style=none] (2) at (0.5, -3) {$X$};
	\end{pgfonlayer}
	\begin{pgfonlayer}{edgelayer}
		\draw (0) to (1.center);
	\end{pgfonlayer}
\end{tikzpicture}
}
  \end{gather*}
\end{depict}

\begin{defn}[Causality] \label{def:causal}
  We say that a morphism is \textit{causal} if it satisfies the following
  condition, where \(\ground\) is the canonical discarding map (supplied by the
  copy-delete category structure) of the appropriate type; compare Observation
  \ref{obs:ground}.
  \[
  \begin{tikzpicture}
	\begin{pgfonlayer}{nodelayer}
		\node [style=none] (0) at (0, 0) {$=$};
		\node [style=arrow box] (1) at (-2, 0) {$c$};
		\node [style=discarder] (2) at (-2, 1.5) {};
		\node [style=discarder] (3) at (2, 1.5) {};
		\node [style=none] (4) at (-2, -1.5) {};
		\node [style=none] (5) at (2, -1.5) {};
	\end{pgfonlayer}
	\begin{pgfonlayer}{edgelayer}
		\draw (4.center) to (1);
		\draw (1) to (2);
		\draw (5.center) to (3);
	\end{pgfonlayer}
\end{tikzpicture}
  \]
\end{defn}

\begin{rmk}
  Observe that, if the monoidal unit is terminal, then every morphism is causal.
\end{rmk}

\begin{depict}[Monoidal symmetry]
  The symmetry of the monoidal structure \(\mathsf{swap}_{XY} : X \otimes Y
  \xto{\sim} Y \otimes X\) is depicted as the swapping of wires, and satisfies
  the equations below. The left says that swapping is an isomorphism; the right
  says that it commutes with copying, making every object a \textit{commutative}
  comonoid:
  
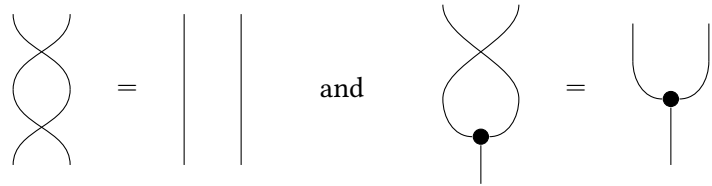
\begin{equation*} %
    \begin{tikzpicture}
	\begin{pgfonlayer}{nodelayer}
		\node [style=none] (0) at (-3, -2) {};
		\node [style=none] (1) at (-3, 0) {};
		\node [style=none] (2) at (-3, 2) {};
		\node [style=none] (3) at (-1.5, 2) {};
		\node [style=none] (4) at (-1.5, 0) {};
		\node [style=none] (5) at (-1.5, -2) {};
		\node [style=none] (6) at (0, 0) {$=$};
		\node [style=none] (7) at (1.5, -2) {};
		\node [style=none] (8) at (1.5, 2) {};
		\node [style=none] (9) at (3, 2) {};
		\node [style=none] (10) at (3, -2) {};
	\end{pgfonlayer}
	\begin{pgfonlayer}{edgelayer}
		\draw [in=270, out=90] (0.center) to (4.center);
		\draw [in=270, out=90] (5.center) to (1.center);
		\draw [in=-90, out=90] (1.center) to (3.center);
		\draw [in=-90, out=90] (4.center) to (2.center);
		\draw (7.center) to (8.center);
		\draw (10.center) to (9.center);
	\end{pgfonlayer}
\end{tikzpicture}
    \hspace{0.06\linewidth}
    \text{and}
    \hspace{0.06\linewidth}
    \begin{tikzpicture}
	\begin{pgfonlayer}{nodelayer}
		\node [style=none] (0) at (0, 0) {$=$};
		\node [style=black copier] (1) at (-2.5, -1.25) {};
		\node [style=none] (2) at (-1.5, -0.25) {};
		\node [style=none] (3) at (-3.5, -0.25) {};
		\node [style=none] (5) at (-2.5, -2.5) {};
		\node [style=none] (6) at (-3.5, 2.25) {};
		\node [style=none] (7) at (-1.5, 2.25) {};
		\node [style=black copier] (8) at (2.5, -0.25) {};
		\node [style=none] (9) at (3.5, 0.75) {};
		\node [style=none] (10) at (1.5, 0.75) {};
		\node [style=none] (11) at (2.5, -2) {};
		\node [style=none] (12) at (1.5, 1.75) {};
		\node [style=none] (13) at (3.5, 1.75) {};
	\end{pgfonlayer}
	\begin{pgfonlayer}{edgelayer}
		\draw (5.center) to (1);
		\draw [in=270, out=180] (1) to (3.center);
		\draw [in=270, out=0] (1) to (2.center);
		\draw [in=-90, out=90, looseness=0.75] (3.center) to (7.center);
		\draw [in=-90, out=90, looseness=0.75] (2.center) to (6.center);
		\draw (11.center) to (8);
		\draw [in=270, out=180] (8) to (10.center);
		\draw [in=270, out=0] (8) to (9.center);
		\draw (10.center) to (12.center);
		\draw (9.center) to (13.center);
	\end{pgfonlayer}
\end{tikzpicture}
  \end{equation*}
\end{depict}

\begin{depict}[Comonoid laws]
  The copy-delete structure \((\bcopier, \ground)\) is required to satisfy the
  comonoid laws, depicted below, of unitality (left depiction) and associativity
  (right depiction):
  
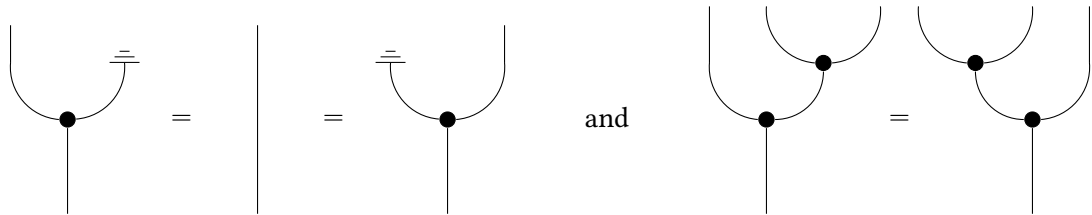
\begin{equation*} %
    \begin{tikzpicture}
	\begin{pgfonlayer}{nodelayer}
		\node [style=none] (0) at (0, 2.5) {};
		\node [style=none] (1) at (0, -2.5) {};
		\node [style=none] (2) at (2, 0) {$=$};
		\node [style=none] (3) at (-2, 0) {$=$};
		\node [style=black copier] (4) at (-5, 0) {};
		\node [style=black copier] (5) at (5, 0) {};
		\node [style=discarder] (6) at (-3.5, 1.5) {};
		\node [style=discarder] (7) at (3.5, 1.5) {};
		\node [style=none] (8) at (6.5, 1.5) {};
		\node [style=none] (9) at (-6.5, 1.5) {};
		\node [style=none] (10) at (-6.5, 2.5) {};
		\node [style=none] (11) at (6.5, 2.5) {};
		\node [style=none] (12) at (-5, -2.5) {};
		\node [style=none] (13) at (5, -2.5) {};
	\end{pgfonlayer}
	\begin{pgfonlayer}{edgelayer}
		\draw (12.center) to (4);
		\draw [bend left=45] (4) to (9.center);
		\draw (9.center) to (10.center);
		\draw [bend right=45] (4) to (6);
		\draw (1.center) to (0.center);
		\draw (13.center) to (5);
		\draw [bend left=45] (5) to (7);
		\draw [bend right=45] (5) to (8.center);
		\draw (8.center) to (11.center);
	\end{pgfonlayer}
\end{tikzpicture}
    \hspace{0.06\linewidth}
    \text{and}
    \hspace{0.06\linewidth}
    \begin{tikzpicture}
	\begin{pgfonlayer}{nodelayer}
		\node [style=black copier] (0) at (-3.5, 0) {};
		\node [style=black copier] (1) at (-2, 1.5) {};
		\node [style=none] (2) at (-5, 1.5) {};
		\node [style=none] (3) at (-5, 3) {};
		\node [style=none] (4) at (-3.5, -2.5) {};
		\node [style=none] (5) at (-3.5, 3) {};
		\node [style=none] (6) at (-0.5, 3) {};
		\node [style=black copier] (7) at (3.5, 0) {};
		\node [style=black copier] (8) at (2, 1.5) {};
		\node [style=none] (9) at (5, 1.5) {};
		\node [style=none] (10) at (5, 3) {};
		\node [style=none] (11) at (3.5, -2.5) {};
		\node [style=none] (12) at (3.5, 3) {};
		\node [style=none] (13) at (0.5, 3) {};
		\node [style=none] (14) at (0, 0) {$=$};
	\end{pgfonlayer}
	\begin{pgfonlayer}{edgelayer}
		\draw (4.center) to (0);
		\draw [bend left=45] (0) to (2.center);
		\draw (2.center) to (3.center);
		\draw [bend right=45] (0) to (1);
		\draw [bend left=45] (1) to (5.center);
		\draw [bend right=45] (1) to (6.center);
		\draw (11.center) to (7);
		\draw [bend right=45] (7) to (9.center);
		\draw (9.center) to (10.center);
		\draw [bend left=45] (7) to (8);
		\draw [bend right=45] (8) to (12.center);
		\draw [bend left=45] (8) to (13.center);
	\end{pgfonlayer}
\end{tikzpicture}
  \end{equation*}
\end{depict}

\begin{depict}[Marginalization of joint states] \label{depic:marginals}
  The discarding maps induce projections \(X \otimes Y \xto{\id \otimes \ground}
  X \otimes I \xto{\sim} X\) and \(X \otimes Y \xto{\ground \otimes \id} I
  \otimes Y \xto{\sim} Y\), with which we can obtain the marginals of joint
  states; compare Definition \ref{def:marginals}. Suppose then that a joint
  state \(\omega : I \klto X \otimes Y\) has marginals \(\omega_1 : I \klto X\)
  and \(\omega_2 : I \klto Y\). Then we have
  
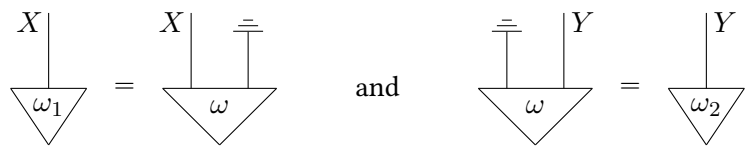
\[
  \begin{tikzpicture}
	\begin{pgfonlayer}{nodelayer}
		\node [style=none] (0) at (1, 0) {};
		\node [style=none] (1) at (2.5, -1.5) {};
		\node [style=none] (2) at (4, 0) {};
		\node [style=none] (3) at (2.5, -0.5) {$\omega$};
		\node [style=none] (4) at (1.75, 0) {};
		\node [style=none] (5) at (3.25, 0) {};
		\node [style=discarder] (6) at (3.25, 1.5) {};
		\node [style=none] (7) at (1.75, 2) {};
		\node [style=none] (8) at (0, 0) {$=$};
		\node [style=none] (9) at (-2, 2) {};
		\node [style=none] (10) at (-2, -0.5) {$\omega_1$};
		\node [style=none] (11) at (-2, 0) {};
		\node [style=none] (12) at (-3, 0) {};
		\node [style=none] (13) at (-2, -1.5) {};
		\node [style=none] (14) at (-1, 0) {};
		\node [style=none] (15) at (-2.5, 1.75) {$X$};
		\node [style=none] (16) at (1.25, 1.75) {$X$};
	\end{pgfonlayer}
	\begin{pgfonlayer}{edgelayer}
		\draw (1.center) to (0.center);
		\draw (1.center) to (2.center);
		\draw (2.center) to (0.center);
		\draw (5.center) to (6);
		\draw (4.center) to (7.center);
		\draw (13.center) to (12.center);
		\draw (13.center) to (14.center);
		\draw (14.center) to (12.center);
		\draw (11.center) to (9.center);
	\end{pgfonlayer}
\end{tikzpicture}
  \hspace{0.06\linewidth}
  \text{and}
  \hspace{0.06\linewidth}
  \begin{tikzpicture}
	\begin{pgfonlayer}{nodelayer}
		\node [style=none] (0) at (-1, 0) {};
		\node [style=none] (1) at (-2.5, -1.5) {};
		\node [style=none] (2) at (-4, 0) {};
		\node [style=none] (3) at (-2.5, -0.5) {$\omega$};
		\node [style=none] (4) at (-1.75, 0) {};
		\node [style=none] (5) at (-3.25, 0) {};
		\node [style=discarder] (6) at (-3.25, 1.5) {};
		\node [style=none] (7) at (-1.75, 2) {};
		\node [style=none] (8) at (0, 0) {$=$};
		\node [style=none] (9) at (2, 2) {};
		\node [style=none] (10) at (2, -0.5) {$\omega_2$};
		\node [style=none] (11) at (2, 0) {};
		\node [style=none] (12) at (3, 0) {};
		\node [style=none] (13) at (2, -1.5) {};
		\node [style=none] (14) at (1, 0) {};
		\node [style=none] (15) at (2.5, 1.75) {$Y$};
		\node [style=none] (16) at (-1.25, 1.75) {$Y$};
	\end{pgfonlayer}
	\begin{pgfonlayer}{edgelayer}
		\draw (1.center) to (0.center);
		\draw (1.center) to (2.center);
		\draw (2.center) to (0.center);
		\draw (5.center) to (6);
		\draw (4.center) to (7.center);
		\draw (13.center) to (12.center);
		\draw (13.center) to (14.center);
		\draw (14.center) to (12.center);
		\draw (11.center) to (9.center);
	\end{pgfonlayer}
\end{tikzpicture}
  \, .
  \]
\end{depict}

\begin{depict}[Generative models]
  A generative model \((\pi, c) : X \klto Y\) induces a joint state \(\omega\)
  on \(X \otimes Y\) by
  
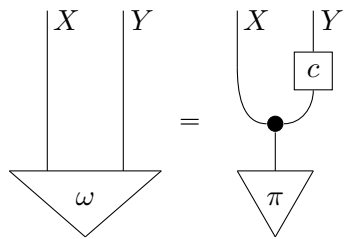
\[
  \begin{tikzpicture}
	\begin{pgfonlayer}{nodelayer}
		\node [style=black copier] (1) at (2.25, 0) {};
		\node [style=none] (2) at (1.25, 3) {};
		\node [style=arrow box] (3) at (3.25, 1.5) {$c$};
		\node [style=none] (4) at (3.25, 3) {};
		\node [style=none] (5) at (1.25, 1.5) {};
		\node [style=none] (6) at (2.25, -2) {$\pi$};
		\node [style=none] (7) at (1.75, 2.75) {$X$};
		\node [style=none] (8) at (3.75, 2.75) {$Y$};
		\node [style=none] (9) at (-4.75, -1.25) {};
		\node [style=none] (10) at (-2.75, -3) {};
		\node [style=none] (11) at (-0.75, -1.25) {};
		\node [style=none] (12) at (-2.75, -2) {$\omega$};
		\node [style=none] (13) at (-3.75, -1.25) {};
		\node [style=none] (14) at (-1.75, -1.25) {};
		\node [style=none] (16) at (-3.75, 3) {};
		\node [style=none] (17) at (-3.25, 2.75) {$X$};
		\node [style=none] (18) at (1.25, -1.25) {};
		\node [style=none] (19) at (3.25, -1.25) {};
		\node [style=none] (20) at (2.25, -3) {};
		\node [style=none] (21) at (2.25, -1.25) {};
		\node [style=none] (22) at (-1.75, 3) {};
		\node [style=none] (23) at (-1.25, 2.75) {$Y$};
		\node [style=none] (24) at (0, 0) {$=$};
	\end{pgfonlayer}
	\begin{pgfonlayer}{edgelayer}
		\draw [in=270, out=-180] (1) to (5.center);
		\draw (5.center) to (2.center);
		\draw (3) to (4.center);
		\draw [in=-90, out=0] (1) to (3);
		\draw (10.center) to (9.center);
		\draw (10.center) to (11.center);
		\draw (11.center) to (9.center);
		\draw (13.center) to (16.center);
		\draw (21.center) to (1);
		\draw (20.center) to (18.center);
		\draw (18.center) to (19.center);
		\draw (19.center) to (20.center);
		\draw (14.center) to (22.center);
	\end{pgfonlayer}
\end{tikzpicture}
  \]
  with marginals \(\pi\) and \(c \klcirc \pi\) given by
  \[
  \begin{tikzpicture}
	\begin{pgfonlayer}{nodelayer}
		\node [style=black copier] (1) at (2, 0) {};
		\node [style=none] (2) at (1, 2.5) {};
		\node [style=arrow box] (3) at (3, 1) {$c$};
		\node [style=discarder] (4) at (3, 2) {};
		\node [style=none] (5) at (1, 1) {};
		\node [style=none] (6) at (2, -2) {$\pi$};
		\node [style=none] (7) at (1.5, 2.25) {$X$};
		\node [style=none] (9) at (-4.5, -1.25) {};
		\node [style=none] (10) at (-2.5, -3) {};
		\node [style=none] (11) at (-0.5, -1.25) {};
		\node [style=none] (12) at (-2.5, -2) {$\omega$};
		\node [style=none] (13) at (-3.5, -1.25) {};
		\node [style=none] (14) at (-1.5, -1.25) {};
		\node [style=none] (16) at (-3.5, 2.5) {};
		\node [style=none] (17) at (-3, 2.25) {$X$};
		\node [style=none] (18) at (1, -1.25) {};
		\node [style=none] (19) at (3, -1.25) {};
		\node [style=none] (20) at (2, -3) {};
		\node [style=none] (21) at (2, -1.25) {};
		\node [style=discarder] (22) at (-1.5, 2) {};
		\node [style=none] (24) at (0, 0) {$=$};
		\node [style=none] (25) at (4.25, 0) {$=$};
		\node [style=none] (26) at (5.5, 2.5) {};
		\node [style=none] (27) at (5.5, -1.25) {};
		\node [style=none] (28) at (4.5, -1.25) {};
		\node [style=none] (29) at (6.5, -1.25) {};
		\node [style=none] (30) at (5.5, -3) {};
		\node [style=none] (31) at (5.5, -2) {$\pi$};
		\node [style=none] (32) at (6, 2.25) {$X$};
	\end{pgfonlayer}
	\begin{pgfonlayer}{edgelayer}
		\draw [in=270, out=-180] (1) to (5.center);
		\draw (5.center) to (2.center);
		\draw (3) to (4);
		\draw [in=-90, out=0] (1) to (3);
		\draw (10.center) to (9.center);
		\draw (10.center) to (11.center);
		\draw (11.center) to (9.center);
		\draw (13.center) to (16.center);
		\draw (21.center) to (1);
		\draw (20.center) to (18.center);
		\draw (18.center) to (19.center);
		\draw (19.center) to (20.center);
		\draw (14.center) to (22);
		\draw (30.center) to (28.center);
		\draw (30.center) to (29.center);
		\draw (29.center) to (28.center);
		\draw (27.center) to (26.center);
	\end{pgfonlayer}
\end{tikzpicture}
  \hspace{0.06\linewidth}
  \text{and}
  \hspace{0.06\linewidth}
  \begin{tikzpicture}
	\begin{pgfonlayer}{nodelayer}
		\node [style=black copier] (1) at (2, 0) {};
		\node [style=discarder] (2) at (1, 2) {};
		\node [style=arrow box] (3) at (3, 1) {$c$};
		\node [style=none] (4) at (3, 2.5) {};
		\node [style=none] (5) at (1, 1) {};
		\node [style=none] (6) at (2, -2) {$\pi$};
		\node [style=none] (7) at (3.5, 2.25) {$Y$};
		\node [style=none] (9) at (-4.5, -1.25) {};
		\node [style=none] (10) at (-2.5, -3) {};
		\node [style=none] (11) at (-0.5, -1.25) {};
		\node [style=none] (12) at (-2.5, -2) {$\omega$};
		\node [style=none] (13) at (-3.5, -1.25) {};
		\node [style=none] (14) at (-1.5, -1.25) {};
		\node [style=discarder] (16) at (-3.5, 2) {};
		\node [style=none] (17) at (-1, 2.25) {$Y$};
		\node [style=none] (18) at (1, -1.25) {};
		\node [style=none] (19) at (3, -1.25) {};
		\node [style=none] (20) at (2, -3) {};
		\node [style=none] (21) at (2, -1.25) {};
		\node [style=none] (22) at (-1.5, 2.5) {};
		\node [style=none] (24) at (0, 0) {$=$};
		\node [style=none] (25) at (4.25, 0) {$=$};
		\node [style=none] (27) at (5.5, -1.25) {};
		\node [style=none] (28) at (4.5, -1.25) {};
		\node [style=none] (29) at (6.5, -1.25) {};
		\node [style=none] (30) at (5.5, -3) {};
		\node [style=none] (31) at (5.5, -2) {$\pi$};
		\node [style=none] (32) at (6, 2.25) {$Y$};
		\node [style=arrow box] (33) at (5.5, 1) {$c$};
		\node [style=none] (34) at (5.5, 2.5) {};
	\end{pgfonlayer}
	\begin{pgfonlayer}{edgelayer}
		\draw [in=270, out=-180] (1) to (5.center);
		\draw (5.center) to (2);
		\draw (3) to (4.center);
		\draw [in=-90, out=0] (1) to (3);
		\draw (10.center) to (9.center);
		\draw (10.center) to (11.center);
		\draw (11.center) to (9.center);
		\draw (13.center) to (16);
		\draw (21.center) to (1);
		\draw (20.center) to (18.center);
		\draw (18.center) to (19.center);
		\draw (19.center) to (20.center);
		\draw (14.center) to (22.center);
		\draw (30.center) to (28.center);
		\draw (30.center) to (29.center);
		\draw (29.center) to (28.center);
		\draw (33) to (34.center);
		\draw (27.center) to (33);
	\end{pgfonlayer}
\end{tikzpicture} \, .
  \]
  (Compare Definition \ref{def:gen-model}.)
\end{depict}

\begin{defn}[Bayesian inversion] \label{def:admit-bayes}
  We say that a channel \(c : X \klto Y\) \textit{admits Bayesian inversion}
  with respect to \(\pi : I \klto X\) if there exists a channel \(c^\dag_\pi : Y
  \klto X\) satisfying the following equation
  \parencite[eq. 5]{Cho2017Disintegration}:
  \begin{equation} \label{eq:bayes-abstr}
    \begin{tikzpicture}
	\begin{pgfonlayer}{nodelayer}
		\node [style=black copier] (1) at (0, 0) {};
		\node [style=none] (2) at (-2, 6) {};
		\node [style=arrow box] (3) at (2, 3.5) {$c$};
		\node [style=none] (4) at (2, 6) {};
		\node [style=none] (5) at (-2, 2) {};
		\node [style=state180] (6) at (0, -4.5) {$\pi$};
		\node [style=none] (7) at (-1.5, 5.75) {$X$};
		\node [style=none] (8) at (2.5, 5.75) {$Y$};
		\node [style=none] (9) at (2, 2) {};
	\end{pgfonlayer}
	\begin{pgfonlayer}{edgelayer}
		\draw [bend left=45] (1) to (5.center);
		\draw (5.center) to (2.center);
		\draw (3) to (4.center);
		\draw (6) to (1);
		\draw [bend right=45] (1) to (9.center);
		\draw (9.center) to (3);
	\end{pgfonlayer}
\end{tikzpicture} \quad = \quad \begin{tikzpicture}
	\begin{pgfonlayer}{nodelayer}
		\node [style=black copier] (0) at (0, 0) {};
		\node [style=none] (1) at (2, 6) {};
		\node [style=arrow box] (2) at (-2, 3.5) {$c^\dag_\pi$};
		\node [style=none] (3) at (-2, 6) {};
		\node [style=none] (4) at (2, 2) {};
		\node [style=state180] (5) at (0, -4.5) {$\pi$};
		\node [style=arrow box] (6) at (0, -2) {$c$};
		\node [style=none] (7) at (-1.5, 5.75) {$X$};
		\node [style=none] (8) at (2.5, 5.75) {$Y$};
		\node [style=none] (9) at (-2, 2) {};
	\end{pgfonlayer}
	\begin{pgfonlayer}{edgelayer}
		\draw [bend right=45] (0) to (4.center);
		\draw (4.center) to (1.center);
		\draw (2) to (3.center);
		\draw (5) to (6);
		\draw (6) to (0);
		\draw [bend left=45] (0) to (9.center);
		\draw (9.center) to (2);
	\end{pgfonlayer}
\end{tikzpicture}
  \end{equation}
  We say that \(c\) admits Bayesian inversion \textit{tout court} if \(c\)
  admits Bayesian inversion with respect to all states \(\pi : I \klto X\) such
  that \(c \klcirc \pi\) has non-empty support. We say that a category
  \(\cat{C}\) admits Bayesian inversion if all its morphisms admit Bayesian
  inversion \textit{tout court}.
\end{defn}

\begin{depict}[Density functions] \label{def:density-functions-graph}
  A channel \(c : X \klto Y\) is said to be \textit{represented by the effect}
  \(p_c : X \otimes Y \klto I\) with respect to \(\mu : I \klto Y\) if
  \[ \begin{tikzpicture}
	\begin{pgfonlayer}{nodelayer}
		\node [style=none] (0) at (-1, -4) {};
		\node [style=arrow box] (1) at (-1, 0) {$c$};
		\node [style=none] (2) at (-1, 4) {};
		\node [style=none] (3) at (-0.5, -3.75) {$X$};
		\node [style=none] (4) at (-0.5, 3.5) {$Y$};
		\node [style=black copier] (5) at (7, 0) {};
		\node [style=none] (6) at (9, 2) {};
		\node [style=state180] (7) at (7, -2) {$\mu$};
		\node [style=none] (8) at (2, 2) {};
		\node [style=none] (9) at (6, 2) {};
		\node [style=none] (10) at (4, 3.5) {};
		\node [style=none] (11) at (4, 2.5) {$p$};
		\node [style=none] (12) at (5, 2) {};
		\node [style=none] (13) at (3, 2) {};
		\node [style=none] (14) at (9, 4) {};
		\node [style=none] (15) at (3, -4) {};
		\node [style=none] (16) at (3.5, -3.75) {$X$};
		\node [style=none] (17) at (9.5, 3.75) {$Y$};
		\node [style=none] (18) at (1, 0) {$=$};
	\end{pgfonlayer}
	\begin{pgfonlayer}{edgelayer}
		\draw (0.center) to (1);
		\draw (1) to (2.center);
		\draw (8.center) to (10.center);
		\draw (10.center) to (9.center);
		\draw (8.center) to (9.center);
		\draw (15.center) to (13.center);
		\draw [bend right=45] (12.center) to (5);
		\draw [bend right=45] (5) to (6.center);
		\draw (5) to (7);
		\draw (6.center) to (14.center);
	\end{pgfonlayer}
\end{tikzpicture}. \]
  We call \(p_c\) a \textit{density function} for \(c\); compare Definition
  \ref{def:density-func}.
\end{depict}

\begin{defn}[Almost-equality] \label{def:almost-eq}
  Given a state \(\pi : I \klto X\), we say that two channels \(c : X \klto Y\)
  and \(d : X \klto Y\) are \(\pi\)\textit{-almost-equal}, denoted \(c
  \overset{\pi}{\sim} d\), if
  \[ \begin{tikzpicture}
	\begin{pgfonlayer}{nodelayer}
		\node [style=black copier] (1) at (0, 0) {};
		\node [style=none] (2) at (-2, 6) {};
		\node [style=arrow box] (3) at (2, 3.5) {$c$};
		\node [style=none] (4) at (2, 6) {};
		\node [style=none] (5) at (-2, 2) {};
		\node [style=state180] (6) at (0, -4.5) {$\pi$};
		\node [style=none] (7) at (-1.5, 5.75) {$X$};
		\node [style=none] (8) at (2.5, 5.75) {$Y$};
		\node [style=none] (9) at (2, 2) {};
	\end{pgfonlayer}
	\begin{pgfonlayer}{edgelayer}
		\draw [bend left=45] (1) to (5.center);
		\draw (5.center) to (2.center);
		\draw (3) to (4.center);
		\draw (6) to (1);
		\draw [bend right=45] (1) to (9.center);
		\draw (9.center) to (3);
	\end{pgfonlayer}
\end{tikzpicture} \ \cong\ \begin{tikzpicture}
	\begin{pgfonlayer}{nodelayer}
		\node [style=black copier] (1) at (0, 0) {};
		\node [style=none] (2) at (-2, 6) {};
		\node [style=arrow box] (3) at (2, 3.5) {$d$};
		\node [style=none] (4) at (2, 6) {};
		\node [style=none] (5) at (-2, 2) {};
		\node [style=state180] (6) at (0, -4.5) {$\pi$};
		\node [style=none] (7) at (-1.5, 5.75) {$X$};
		\node [style=none] (8) at (2.5, 5.75) {$Y$};
		\node [style=none] (9) at (2, 2) {};
	\end{pgfonlayer}
	\begin{pgfonlayer}{edgelayer}
		\draw [bend left=45] (1) to (5.center);
		\draw (5.center) to (2.center);
		\draw (3) to (4.center);
		\draw (6) to (1);
		\draw [bend right=45] (1) to (9.center);
		\draw (9.center) to (3);
	\end{pgfonlayer}
\end{tikzpicture} \, . \]
\end{defn}

\begin{prop}[Composition preserves almost-equality] \label{prop:comp-preserve-almost-eq}
  If \(c \overset{\pi}{\sim} d\), then \(f \klcirc c \overset{\pi}{\sim} f
  \klcirc d\).
  \begin{proof}
    Immediate from the definition of almost-equality.
  \end{proof}
\end{prop}

\begin{prop}[Bayesian inverses are almost-equal] \label{prop:bayes-almost-equal}
  Suppose \(\alpha : Y \klto X\) and \(\beta : Y \klto X\) are both Bayesian
  inversions of the channel \(c : X \klto Y\) with respect to \(\pi : I \klto
  X\). Then \(\alpha \overset{c \klcirc \pi}{\sim} \beta\).
  \begin{proof}
    Immediate from Equation \eqref{eq:bayes-abstr}.
  \end{proof}
\end{prop}

\subsection{Lenses for dependent bidirectional processes} \label{sec:bg-lens}

The Bayesian inversion of a stochastic channel \(c : X \klto Y\) is a family of
channels \(c^\dag_\pi : Y \klto X\) in the opposite direction, indexed by states
on \(X\). Pairs \((c, c^\dag)\) of a morphism \(c\) with a \(c\)-dependent
opposite morphism \(c^\dag\) often fall into the compositional `lens' pattern,
and the Bayesian case is no exception.  In this section, we sketch the basic
theory of lenses, and refer the reader to our preprint
\parencite{Smithe2020Bayesian} for further exposition. The central element of
the structure is a (pseudo)functor picking out, for each morphism \(c\), the
category in which \(c^\dag\) lives. With this piece to hand, an entire
corresponding category of lenses can defined most concisely.

\begin{defn}[{\textcite[Def. 3.3]{Spivak2019Generalized}}]
  The category \(\Cat{GrLens}_F\) of \textit{Grothendieck lenses} for a
  pseudofunctor \(F : \cat{C}\op \to \Cat{Cat}\) is the total category of the
  Grothendieck construction for the pointwise opposite of \(F\).
\end{defn}

\begin{prop}[\(\Cat{GrLens}_F\) is a category] \label{prop:grlens-cat}
  The objects \((\Cat{GrLens}_F)_0\) of \(\Cat{GrLens}_F\) are (dependent) pairs
  \((C, X)\) with \(C : \cat{C}\) and \(X : F(C)\), and its hom-sets
  \(\Cat{GrLens}_F \big( (C, X), (C', X') \big)\) are dependent sums
  \begin{equation*}
    \Cat{GrLens}_F \big( (C, X), (C', X') \big)
    = \sum_{f \, : \, \cat{C}(C, C')} F(C) \big( F(f)(X'), X \big)
  \end{equation*}
  so that a morphism \((C, X) \lensto (C', X')\) is a pair \((f, f^\dag)\) of
  \(f : \cat{C}(C, C')\) and \(f^\dag : F(C) \big( F(f)(X'), X \big)\). We call
  such pairs \textit{Grothendieck lenses for} \(F\) or \(F\)\textit{-lenses}.
  \begin{proof}[Proof sketch]
    The identity Grothendieck lens on \((C, X)\) is \(\id_{(C, X)} = (\id_C,
    \id_X)\). Sequential composition is as follows. Given \((f, f^\dag) : (C, X)
    \lensto (C', X')\) and \((g, g^\dag) : (C', X') \lensto (D, Y)\), their
    composite \((g, g^\dag) \lenscirc (f, f^\dag)\) is defined to be the lens
    \(\big(g \circ f, F(f)(g^\dag) \big) : (C, X) \lensto (D,
    Y)\). Associativity and unitality of composition follow from functoriality
    of \(F\).
  \end{proof}
\end{prop}

\begin{defn} \label{def:simp-lens}
  Suppose \(F(C)_0 = \cat{C}_0\), with \(F : \cat{C}\op \to \Cat{Cat}\) a
  pseudofunctor. Define \(\Cat{SimpGrLens}_F\) to be the full subcategory of
  \(\Cat{GrLens}_F\) whose objects are duplicate pairs \((C, C)\) of objects
  \(C\) in \(\cat{C}\). We call \(\Cat{SimpGrLens}_F\) the category of
  \textit{simple} \(F\)-lenses. More generally, any lens between such duplicate
  pairs will be called a \textit{simple lens}. Since duplicating the objects in
  the pairs \((X,X)\) is redundant, we will write the objects simply as \(X\).
\end{defn}

Another name for a pseudofunctor \(F : \cat{C}\op \to \Cat{Cat}\) is an
\textit{indexed category}. When \(\cat{C}\) is a monoidal category with
which \(F\) is appropriately compatible, then we can `upgrade' the notions of
indexed category and Grothendieck construction accordingly. In this work, the
domain categories \(\cat{C}\) are only trivially bicategories, and the
pseudofunctors \(F\) are really just functors; we will restrict ourselves
to the 1-categorical case of monoidal indexed categories, too.

\begin{defn}[{\textcite[{\S}3.2]{Moeller2018Monoidal}}]
  Suppose $(\cat{C},\otimes,I)$ is a monoidal category.
  We say that $F$ is a \textit{monoidal indexed category} when $F$ is a weak lax monoidal functor
  $(F,\mu,\mu_0) : (\cat{C}\op,\otimes\op,I)\to(\Cat{Cat},\times,\Cat{1})$.
  This means that the laxator $\mu$ is given by a natural family of functors $\mu_{A,B} : FA\times FB \Rightarrow F(A\otimes B)$ along with, for any morphisms $f:A\to A'$ and $g:B\to B'$ in $\cat{C}$, a natural isomorphism $\mu_{f,g} : \mu_{A',B'} \circ \left(Ff\times Fg\right) \Rightarrow F(f\otimes g) \circ \mu_{A,B}$. The laxator and the unitor $\mu_0 : \Cat{1}\to FI$ together satisfy standard axioms of associativity and unitality.
  All told, this structure makes \((F, \otimes, \mu, I, \mu_0)\) into a pseudomonoid in the monoidal 2-category of indexed categories and indexed functors.
\end{defn}

\begin{prop}[{\textcite[{\S}6.1]{Moeller2018Monoidal}}] \label{prop:monoidal-gr}
  Suppose \((F, \mu, \mu_0) : (\cat{C}\op, \otimes\op, I) \to (\Cat{Cat},
  \times, \Cat{1})\) is a monoidal indexed category. Then the total
  category of the Grothendieck construction \(\int F\) obtains a monoidal
  structure \((\otimes_\mu, I_\mu)\). On objects, define
  \[
  (C, X) \otimes_\mu (D, Y) := \big(C \otimes D, \mu_{CD}(X, Y)\big)
  \]
  where \(\mu_{CD} : FC \times FD \to F(C \otimes D)\) is the component of
  \(\mu\) at \((C, D)\). On morphisms \((f, f^\dag) : (C, X) \lensto (C', X')\)
  and \((g, g^\dag) : (D, Y) \lensto (D', Y')\), define
  \[
  (f,f^\dag) \otimes_\mu (g, g^\dag) := \big(f \otimes g, \mu_{CD}(f^\dag, g^\dag)\big) \, .
  \]
  The monoidal unit \(I_\mu\) is defined to be the object \(I_\mu := \big(I,
  \mu_0(\ast)\big)\). Writing \(\lambda : I \otimes (-) \Rightarrow (-)\) and
  \(\rho : C \otimes (-) \Rightarrow (-)\) for the left and right unitors of the
  monoidal structure on \(\cat{C}\), the left and right unitors in \(\int F\)
  are given by \((\lambda, \id)\) and \((\rho, \id)\) respectively. Writing
  \(\alpha\) for the associator of the monoidal structure on \(\cat{C}\), the
  associator in \(\int F\) is given by \((\alpha, \id)\).
\end{prop}

\begin{cor} \label{cor:monoidal-grlens}
  When \(F : \cat{C}\op \to \Cat{Cat}\) is equipped with a (weak) lax monoidal
  structure \((\mu, \mu_0)\), its category of lenses \(\Cat{GrLens}_F\) becomes
  a monoidal category \((\Cat{GrLens}_F, \otimes'_\mu, I_\mu)\). On objects
  \(\otimes'_\mu\) is defined as \(\otimes_\mu\) in Proposition
  \ref{prop:monoidal-gr}, as is \(I_\mu\). On morphisms \((f, f^\dag) : (C, X)
  \lensto (C', X')\) and \((g, g^\dag) : (D, Y) \lensto (D', Y')\), define
  \[
  (f,f^\dag) \otimes'_\mu (g, g^\dag) := \big(f \otimes g, \mu_{CD}^{\mathrm{op}}(f^\dag, g^\dag)\big)
  \]
  where \(\mu_{CD}^{\mathrm{op}} : F(C)\op \times F(D)\op \to F(C \otimes
  D)\op\) is the pointwise opposite of \(\mu_{CD}\). The associator and unitors
  are defined as in Proposition \ref{prop:monoidal-gr}.
\end{cor}

\begin{rmk} \label{rmk:optics}
  An alternative perspective on lenses is given by the family of structures known as \textit{optics} \parencite{Clarke2020Profunctor}, which generalize lenses using the kind of actegorical machinery to which we now turn.
  This machinery puts the categories of forwards and backwards maps on equal footing, unlike the fibrational machinery developed here (which privileges the base category), at the cost of the explicit dependence structure and somewhat heavier categorical tooling.
  The synthesis of Grothendieck lenses and optics, \textit{dependent optics}, has recently been articulated \parencite{Vertechi2022Dependent,Capucci2022Seeing,Braithwaite2021Fibre}; while powerful, this structure demands both the machinery of fibrations and of actegories.
  For its relative simplicity, we therefore stick to the fibrational lens perspective in this paper.
\end{rmk}

\subsection{Categories with parameters} \label{sec:bg-para}

In many applications, we will be interested in \textit{cybernetic} systems,
where a single system might have some freedom in the choice of forwards and
backwards channel; consider the synaptic strengths or weights of a neural
network, which change as the system learns about the world, affecting the
predictions it makes and actions it takes. This freedom is well modelled by
equipping the morphisms of a category with parameters, and gives rise to a
notion of \textit{parameterized category}. In general, the parameterization may
have different structure to the processes at hand, and so we describe the
`actegorical' situation in which a category of parameters \(\cat{M}\)
\textit{acts on} on a category of processes \(\cat{C}\), generating a category
of parameterized processes.

\begin{defn}[\(\cat{M}\)-actegory]
  Suppose \(\cat{M}\) is a monoidal category with tensor \(\boxtimes\) and unit
  object \(I\). We say that \(\cat{C}\) is a \textit{left}
  \(\cat{M}\)\textit{-actegory} when \(\cat{C}\) is equipped with a functor
  \(\odot : \cat{M} \to \Cat{Cat}(\cat{C}, \cat{C})\) called the \textit{action}
  along with natural unitor and associator isomorphisms \(\lambda^\odot_X : I
  \odot X \xto{\sim} X\) and \(a^\odot_{M,N,X} : (M \boxtimes N) \odot X
  \xto{\sim} M \odot (N \odot X)\) compatible with the monoidal structure of
  \((\cat{M}, \boxtimes, I)\).
\end{defn}

\begin{prop}[\textcite{Capucci2021Towards}] \label{prop:para-bicat}
  Let \((\cat{C}, \odot, \lambda^\odot, a^\odot)\) be an \((\cat{M}, \boxtimes,
  I)\)-actegory. Then there is a bicategory of \(\cat{M}\)-parameterized
  morphisms in \(\cat{C}\), denoted \(\para(\odot)\). Its objects are those of
  \(\cat{C}\). For each pair of objects \(X,Y\), the set of 1-cells is defined
  as \(\para(\odot)(X, Y) := \sum_{M:\cat{M}} \cat{C}(M \odot X, Y)\); we denote
  an element \((M, f)\) of this set by \(f : X \xto{M} Y\). Given 1-cells \(f :
  X \xto{M} Y\) and \(g : Y \xto{N} Z\), their composite \(g \circ f : X \xto{N
    \boxtimes M} Z\) is the following morphism in \(\cat{C}\):
  \[
  (N \boxtimes M) \odot X \xto{a^\odot_{N,M,X}} N \odot (M \odot X) \xto{\id_N \odot f} N \odot Y \xto{g} Z
  \]
  Given 1-cells \(f : X \xto{M} Y\) and \(f' : X \xto{M'} Y\), a 2-cell \(\alpha
  : f \Rightarrow f'\) is a morphism \(\alpha : M \to M'\) in \(\cat{M}\) such
  that \(f = f' \circ (\alpha \odot \id_X)\) in \(\cat{C}\); identities and
  composition of 2-cells are as in \(\cat{C}\).
\end{prop}

\begin{prop} \label{prop:para-tensor}
  When \(\cat{C}\) is equipped with both a symmetric monoidal structure
  \((\otimes, I)\) and an \((\cat{M}, \boxtimes, I)\)-actegory structure, and
  there is a natural isomorphism \(\mu^r_{M,X,Y} : M \odot (X \otimes Y)
  \xto{\sim} X \otimes (M \odot Y)\) called the \textit{(right)
    costrength}\footnote{So named for its similarity to the (co)strength of
  certain (co)monads; and similarly, too, our costrength should satisfy certain
  standard coherence laws which we omit here.}, the symmetric monoidal structure
  \((\otimes, I)\) lifts to \(\para(\odot)\). First, from the symmetry of
  \(\otimes\), one obtains a (left) costrength, \(\mu^l_{M,X,Y} : M \odot (X
  \otimes Y) \xto{\sim} (M \odot X) \otimes Y\). Second, using the two
  costrengths and the associator of the actegory structure, one obtains a
  natural isomorphism \(\iota_{M,N,X,Y} : (M \boxtimes N) \odot (X \otimes Y)
  \xto{\sim} (M \odot X) \otimes (N \odot Y)\) called the
  \textit{interchanger}. The tensor of objects in \(\para(\odot)\) is then
  defined as the tensor of objects in \(\cat{C}\), and the tensor of morphisms
  (1-cells) \(f : X \xto{M} Y\) and \(g : A \xto{N} B\) is given by the
  composite
  \[
  f \otimes g : X \otimes A \xto{M \boxtimes N} Y \otimes B
  \; := \;
  (M \boxtimes N) \odot (X \otimes A) \xto{\iota_{M,N,X,A}}
  (M \odot A) \otimes (N \odot A) \xto{f \otimes g} Y \otimes B \, .
  \]
\end{prop}

In many simple cases, the parameters will live in the same category as the
morphisms being parameterized; this is formalized by the following proposition.

\begin{prop} \label{prop:para-self}
  If \((\cat{C}, \otimes, I)\) is a monoidal category, then it induces a
  parameterization \(\para(\otimes)\) on itself. For each \(M,X,Y : \cat{C}\),
  the morphisms \(X \xto{M} Y\) of \(\para(\otimes)\) are the morphisms \(M
  \otimes X \to Y\) in \(\cat{C}\).
\end{prop}

\begin{notation}
  When considering the self-paramterization induced by a monoidal category
  \((\cat{C}, \otimes, I)\), we will often write \(\para(\cat{C})\) instead of
  \(\para(\otimes)\).
\end{notation}

It will frequently be the case that we do not in fact need the whole bicategory
structure. The following proposition tells us that we can also just work
1-categorically, as long as we work with equivalence classes of
isomorphically-parameterized maps, in order that composition is suffiently
strictly associative.

\begin{prop} \label{prop:para-1cat}
  Each bicategory \(\para(\odot)\) induces a 1-category \(\para(\odot)_1\) by
  forgetting the bicategorical structure. The hom sets \(\para(\odot)_1(X, Y)\)
  are given by \(U \para(\odot)(X, Y) / \sim\) where \(U\) is the forgetful
  functor \(U : \Cat{Cat} \to \Set\) and \(f \sim g\) if and only if there is
  some 2-cell \(\alpha : f \Rightarrow g\) that is an isomorphism. We call
  \(\para(\odot)_1\) the \textit{1-categorical truncation} of
  \(\para(\odot)\). When \(\para(\odot)\) is monoidal, so is \(\para(\odot)_1\).
\end{prop}

\section{The bidirectional structure of Bayesian updating} \label{sec:buco}

In this section, we define a collection of indexed categories, each denoted
\(\mathsf{Stat}\), whose morphisms can be seen as ``generalized Bayesian
inversions''. Following Proposition \ref{prop:grlens-cat}, these induce
corresponding categories of lenses which we call \textit{Bayesian lenses}. We
show that, for the subcategories of \textit{exact} Bayesian lenses whose
backward channels correspond to `exact' Bayesian inversions, the Bayesian
inversion of a composite of forward channels is given (up to almost-equality) by
the lens composite of the corresponding backward channels. This justifies
calling these lenses `Bayesian', and provides the foundation for the study of
approximate (non-exact) Bayesian inversion in the subsequent section.

\begin{rmk*}
  An alternative account of Bayesian lenses, from an `optical' perspective, is
  told in the preprint \parencite{Smithe2020Bayesian}.
\end{rmk*}

\subsection{State-dependent channels}

A channel \(c : X \klto Y\) admitting a Bayesian inversion induces a family of
inverse channels \(c^\dag_\pi : Y \klto X\), indexed by `prior' states \(\pi : 1
\klto X\). Making the state-dependence explicit, in typical cases where \(c\) is
a probability kernel we obtain a measurable function \(c^\dag : \Giry X \times Y
\to \Giry X\). In more general situations, we obtain a morphism \(c^\dag :
\cat{C}(I, X) \to \cat{C}(Y, X)\) in the base of enrichment of the monoidal
category \((\cat{C}, \otimes, I)\) of \(c\). We call morphisms of this general
type \textit{state-dependent channels}, and structure the indexing as an indexed
category.

\begin{defn} \label{def:stat-cat}
  Let \((\cat{C}, \otimes, I)\) be a monoidal category enriched in a Cartesian
  category \(\Cat{V}\). Define the \(\cat{C}\)\textit{-state-indexed} category
  \(\Fun{Stat}: \cat{C}\op \to \Cat{V\mdash Cat}\) as follows.
  \begin{align}
    \Fun{Stat} \;\; : \;\; \cat{C}\op \; & \to \; \Cat{V\mdash Cat} \nonumber \\
    X & \mapsto \Fun{Stat}(X) := \quad \begin{pmatrix*}[l]
      & \Fun{Stat}(X)_0 & := \quad \;\;\; \cat{C}_0 \\
      & \Fun{Stat}(X)(A, B) & := \quad \;\;\; \Cat{V}(\cat{C}(I, X), \cat{C}(A, B)) \\
      \id_A \: : & \Fun{Stat}(X)(A, A) & := \quad
      \left\{ \begin{aligned}
        \id_A : & \; \cat{C}(I, X)     \to     \cat{C}(A, A) \\
        & \quad\;\;\: \rho \quad \mapsto \quad \id_A
      \end{aligned} \right. \label{eq:stat} \\
    \end{pmatrix*} \\ \nonumber \\
    f : \cat{C}(Y, X) & \mapsto \begin{pmatrix*}[c]
      \Fun{Stat}(f) \; : & \Fun{Stat}(X) & \to & \Fun{Stat}(Y) \vspace*{0.5em} \\
      & \Fun{Stat}(X)_0 & = & \Fun{Stat}(Y)_0 \vspace*{0.5em} \\
      & \Cat{V}(\cat{C}(I, X), \cat{C}(A, B)) & \to & \Cat{V}(\cat{C}(I, Y), \cat{C}(A, B)) \vspace*{0.125em} \\
      & \alpha & \mapsto & f^\ast \alpha : \big( \, \sigma : \cat{C}(I, Y) \, \big) \mapsto \big( \, \alpha(f \klcirc \sigma) : \cat{C}(A, B) \, \big)
    \end{pmatrix*} \nonumber
  \end{align}
  Composition in each fibre \(\Fun{Stat}(X)\) is as in \(\cat{C}\). Explicitly,
  indicating morphisms \(\cat{C}(I, X) \to \cat{C}(A, B)\) in \(\Fun{Stat}(X)\)
  by \(A \xklto{X} B\), and given \(\alpha : A \xklto{X} B\) and \(\beta : B
  \xklto{X} C\), their composite is \(\beta \circ \alpha : A \xklto{X} C := \rho
  \mapsto \beta(\rho) \klcirc \alpha(\rho)\), where here we indicate composition
  in \(\cat{C}\) by \(\klcirc\) and composition in the fibres \(\Fun{Stat}(X)\)
  by \(\circ\). Given \(f : Y \klto X\) in \(\cat{C}\), the induced functor
  \(\Fun{Stat}(f) : \Fun{Stat}(X) \to \Fun{Stat}(Y)\) acts by pullback.
\end{defn}

\begin{notation}
  Just as we wrote \(X \xto{M} Y\) for an \(M\)-parameterized morphism in
  \(\cat{C}(M \odot X, Y)\) (see Proposition \ref{prop:para-bicat}), we write
  \(A \xklto{X} B\) for an \(X\)-state-dependent morphism in
  \(\Cat{V}\big(\cat{C}(I, X), \cat{C}(A, B)\big)\). Given a state \(\rho\) in
  \(\cat{C}(I, X)\) and an \(X\)-state-dependent morphism \(f : A \xklto{X} B\),
  we write \(f_\rho\) for the resulting morphism in \(\cat{C}(A,B)\).
\end{notation}

\begin{rmk} \label{rmk:stat-para-prox}
  We can thus think of \(\Fun{Stat}\) as a kind of `external' parameterization
  of channels in \(\cat{C}\), here by states in \(\cat{C}\). Generalizing this
  notion to ``external parameterization by \(\Cat{V}\) objects'' gives rise to a
  a cousin of the \(\para\) construction: a notion of \textit{category by
    proxy}, denoted \(\Cat{Prox}\), of which \(\Fun{Stat}\) is a (fibred)
  special case. Correspondingly, we think of \(\para\) as capturing
  parameterization `internal' to \(\cat{C}\). Consider for instance the case
  where \(\cat{M}\) acts on \(\cat{C}\) by \(\odot : \cat{M} \times \cat{C} \to
  \cat{C}\) and suppose that, for all \(A : \cat{C}\), the functors \((-)\odot A
  : \cat{M} \to \cat{C}\) are left adjoint to \([A,-] : \cat{C} \to
  \cat{M}\). Then we have \(\para(\odot)(A, B) \cong \sum_{M:\cat{M}}
  \cat{M}\big(M,[A,B]\big)\), which is strongly reminiscent of the definition of
  \(\Fun{Stat}\), and comes close to a definition of \(\Cat{Prox}\). The
  connections between these constructions are a matter of on-going study; a
  summary is available in \textcite{Capucci2021Parameterized}.
\end{rmk}

\begin{prop} \label{prop:stat-lax}
  \(\Fun{Stat}\) is lax monoidal. The components \(\mu_{XY} : \Fun{Stat}(X)
  \times \Fun{Stat}(Y) \to \Fun{Stat}(X \otimes Y)\) of the laxator are defined
  on objects by \(\mu_{XY}(A, A') := A \otimes A'\) and on morphisms \(f : A
  \xklto{X} B\) and \(f' : A' \xklto{Y} B'\) by \(\mu_{XY}(f,f') := f \otimes f'
  : A \otimes A' \xklto{X \otimes Y} B \otimes B'\), where \(f \otimes f'\) is
  the \(\Cat{V}\)-morphism \(\cat{C}(I, X\otimes Y) \to \cat{C}(A\otimes A',
  B\otimes B') : \omega \mapsto f_{\omega_X} \otimes f'_{\omega_Y}\). Here,
  \(\omega_X\) and \(\omega_Y\) are the \(X\) and \(Y\) marginals of \(\omega\),
  defined by \(\omega_X := \mathsf{proj}_X \klcirc \omega\) and \(\omega_Y :=
  \mathsf{proj}_Y \klcirc \omega\) (see Definition \ref{def:marginals} and
  Depiction \ref{depic:marginals}). The unit \(\mu_0 : \Cat{1} \to
  \Fun{Stat}(I)\) of the lax monoidal structure is the functor mapping the
  unique object \(1 : \Cat{1}\) to \(I : \Fun{Stat}(I)\).
\end{prop}

\begin{ex}
  The category \(\Cat{Meas}\) of general measurable spaces is not Cartesian
  closed, as there is no general way to make the evaluation maps \(\Cat{Meas}(X,
  Y) \times X \to Y\) measurable, meaning that if we take \(\cat{C} =
  \Kl(\Giry)\) above, then we are forced to take \(\Cat{V} = \Set\). In turn,
  this makes the inversion maps \(c^\dag : \Kl(\Giry)(1, X) \to
  \Kl(\Giry)(Y,X)\) into mere functions. We can salvage measurability by working
  instead with \(\Kl(\Qa)\), where \(\Qa : \Cat{QBS} \to \Cat{QBS}\) is the
  analogue of the Giry monad for quasi-Borel spaces
  \parencite{Heunen2017Convenient}. The category \(\Cat{QBS}\) is indeed
  Cartesian closed, and \(\Kl(\Qa)\) is enriched in \(\Cat{QBS}\), so that we
  can instantiate \(\Fun{Stat}\) there, and the corresponding inversion maps are
  accordingly measurable. Moreover, there is a quasi-Borel analogue of the
  notion of s-finite kernel \parencite[{\S 11}]{Vakar2018S}, with which we can
  define a variant of the category \(\Cat{sfKrn}\).
\end{ex}

\begin{rmk}
  When \(\cat{C}\) is a Kleisli category \(\Kl(T)\), it is of course possible to
  define a variant of \(\Fun{Stat}\) on the other side of the
  product-exponential adjunction, with state-dependent morphisms \(A \xklto{X}
  B\) having the types \(TX \times A \to TB\). This avoids the technical
  difficulties sketched in the preceding example at the cost of requiring a
  monad \(T\). However, the exponential form makes for better exegesis, and so
  we will stick to that.
\end{rmk}

\subsection{Bayesian lenses} \label{sec:buco-bayes-lens}

We define the category of Bayesian lenses in \(\cat{C}\) to be the category of
Grothendieck \(\Fun{Stat}\)-lenses.

\begin{defn} \label{def:stat-lens}
  The category \(\BLens{C}\) of Bayesian lenses in \(\cat{C}\) is the category
  \(\Cat{GrLens}_{\Fun{Stat}}\) of Grothendieck lenses for the functor
  \(\Fun{Stat}\). A \textit{Bayesian lens} is a morphism in \(\BLens{C}\). Where
  the category \(\cat{C}\) is evident from the context, we will just write
  \(\Cat{BayesLens}\).
\end{defn}

Unpacking this definition, we find that the objects of \(\BLens{C}\) are pairs
\((X, A)\) of objects of \(\cat{C}\). Morphisms (that is, Bayesian lenses) \((X,
A) \lensto (Y, B)\) are pairs \((c, c^\dag)\) of a channel \(c : X \klto Y\) and
a ``generalized Bayesian inversion'' \(c^\dag : B \xklto{X} A\); that is,
elements of the hom objects
\begin{align*}
  \BLens{C}\big((X,A),(Y,B)\big)
  :&= \Cat{GrLens}_\Fun{Stat} \big((X,A),(Y,B)\big) \\
  &\cong \cat{C}(X, Y) \times \Cat{V} \big( \cat{C}(I, X), \cat{C}(B, A) \big) \, .
\end{align*}
The identity Bayesian lens on \((X, A)\) is \((\id_X, \id_A)\), where by abuse
of notation \(\id_A : \cat{C}(I, Y) \to \cat{C}(A, A)\) is the constant map
\(\id_A\) defined in Equation \eqref{eq:stat} that takes any state on \(Y\) to
the identity on \(A\).

The sequential composite \((d, d^\dag) \lenscirc (c, c^\dag)\) of \((c, c^\dag)
: (X, A) \lensto (Y, B)\) and \((d, d^\dag) : (Y, B) \lensto (Z, C)\) is the
Bayesian lens \(\big( (d \klcirc c), (c^\dag \circ c^\ast d^\dag) \big) : (X, A)
\lensto (Z, C)\) where \((c^\dag \circ c^\ast d^\dag) : C \xklto{X} A\) takes a
state \(\pi : I \klto X\) to the channel \(c^\dag_{\pi} \klcirc \d^\dag_{c
  \klcirc \pi} : C \klto A\).

\begin{defn}
  Given a Bayesian lens \((c,c') : (X,A)\lensto (Y,B)\), we will call \(c\) its
  \textit{forwards} or \textit{prediction} channel and \(c'\) its
  \textit{backwards} or \textit{update} channel (even though \(c'\) is really a
  family of channels).
\end{defn}

\begin{rmk}
  Note that the definition of \(\Fun{Stat}\) and hence the definition of
  \(\BLens{C}\) do not require \(\cat{C}\) to be a copy-delete category, even
  though our motivating categories of stochastic channels are; all that is
  required for the definition is that \(\cat{C}\) is monoidal.
\end{rmk}

\begin{rmk} \label{rmk:cartesian-lens}
  On the other hand, the structure of \(\cat{C}\) might be stronger than merely
  monoidal. For instance, when \(\cat{C}\) is Cartesian closed, then we can take
  \(\Cat{V} = \cat{C}\) and the monoidal structure to be the categorical product
  \((\times, 1)\). Then a Bayesian lens \((X,A)\lensto(Y,B)\) is equivalently
  given by a pair of a forwards map \(X\to Y\) and a backwards map \(X\times B
  \to A\). We call such lenses \textit{Cartesian}, and they characterize the
  original `lens' notion; see Remark \ref{rmk:lens-laws}.
\end{rmk}

\begin{prop} \label{prop:blens-monoidal}
  \(\BLens{C}\) is a monoidal category, with structure \(\big((\otimes,
  (I,I)\big)\) inherited from \(\cat{C}\). On objects, define \((A, A') \otimes
  (B, B') := (A \otimes A', B \otimes B')\). On morphisms \((f, f^\dag) : (X, A)
  \lensto (Y, B)\) and \((g, g^\dag) : (X', A') \lensto (Y', B')\), define \((f,
  f^\dag) \otimes (g, g^\dag) := (f \otimes g, f^\dag \otimes g^\dag)\), where
  \(f^\dag \otimes g^\dag : B \otimes B' \xklto{X \otimes X'} A \otimes A'\)
  acts on states \(\omega : I \klto X \otimes X'\) to return the channel
  \(f^\dag_{\omega_X} \otimes g^\dag_{\omega_X'}\), following the definition of
  the laxator \(\mu\) in Proposition \ref{prop:stat-lax}. The monoidal unit in
  \(\BLens{C}\) is the pair \((I,I)\) duplicating the unit in \(\cat{C}\). When
  \(\cat{C}\) is moreover symmetric monoidal, so is \(\BLens{C}\).
  \begin{proof}[Proof sketch]
    The main result is immediate from Proposition \ref{prop:stat-lax} and
    Corollary \ref{cor:monoidal-grlens}. When \(\otimes\) is symmetric in
    \(\cat{C}\), the symmetry lifts to the fibres of \(\Fun{Stat}\) and hence to
    \(\BLens{C}\).
  \end{proof}
\end{prop}

\begin{rmk}
  Although \(\BLens{C}\) is a monoidal category, it does not inherit a copy-delete structure from \(\cat{C}\), owing to the bidirectionality of its component morphisms.
  To see this, we can consider morphisms into the monoidal unit \((I,I)\), and find that there is generally no canonical discarding map.
  For instance, a morphism \((X,A)\lensto(I,I)\) consists in a pair of a channel \(X\klto I\) (which may indeed be a discarding map) and a state-dependent channel \(I \xklto{X} A\), for which there is generally no suitable choice satisfying the comonoid laws.
  Note, however, that a lens of the type \((X,I)\lensto(I,B)\) might indeed act by discarding, since we can choose the constant state-dependent channel \(B\xklto{X} I\) on the discarding map \(\ground : B\klto I\).
  By contrast, the Grothendieck category \(\int \Fun{Stat}\) \textit{is} a copy-delete category, as the morphisms \((X,A)\to(I,I)\) in \(\int\Fun{Stat}\) are pairs \(X\klto I\) and \(A\xklto{X} I\), and so for both components we can choose morphisms witnessing the comonoid structure.
\end{rmk}

\subsection{Bayesian updates compose optically}

In this section we prove the fundamental result on which the development of
statistical games rests: that the inversion of a composite channel is given up
to almost-equality by the lens composite of the backwards components of the
associated `exact' Bayesian lenses.

\begin{defn}
  Let \((c, c^\dag) : (X, X) \lensto (Y, Y)\) be a Bayesian lens. We say that
  \((c, c^\dag)\) is \textit{exact} if \(c\) admits Bayesian inversion and, for
  each \(\pi : I \klto X\) such that \(c \klcirc \pi\) has non-empty support,
  \(c\) and \(c^\dag_\pi\) together satisfy equation \eqref{eq:bayes-abstr}.
  Bayesian lenses that are not exact are said to be \textit{approximate}.
\end{defn}

\begin{thm} \label{thm:buco}
  Let \((c, c^\dag)\) and \((d, d^\dag)\) be sequentially composable exact
  Bayesian lenses. Then the contravariant component of the composite lens \((d,
  d^\dag) \lenscirc (c, c^\dag) = (d \klcirc c, c^\dag \circ c^\ast d^\dag)\)
  is, up to \(d \klcirc c \klcirc \pi\)-almost-equality, the Bayesian inversion
  of \(d \klcirc c\) with respect to any state \(\pi\) on the domain of \(c\)
  such that \(c \klcirc \pi\) has non-empty support. That is to say,
  \emph{Bayesian updates compose optically}: \((d \klcirc c)^\dag_\pi \overset{d
    \klcirc c \klcirc \pi}{\sim} c^\dag_\pi \klcirc d^\dag_{c \klcirc \pi}\).

  \begin{proof}
    Suppose \(c^\dag_\pi : Y \klto X\) is the Bayesian inverse of \(c : X \klto
    Y\) with respect to \(\pi : I \klto X\). Suppose also that \(d^\dag_{c
      \klcirc \pi} : Z \klto Y\) is the Bayesian inverse of \(d : Y \klto X\)
    with respect to \(c \klcirc \pi : I \klto Y\), and that \((d \klcirc
    c)^\dag_\pi : Z \klto X\) is the Bayesian inverse of \(d \klcirc c : X \klto
    Z\) with respect to \(\pi : I \klto X\):
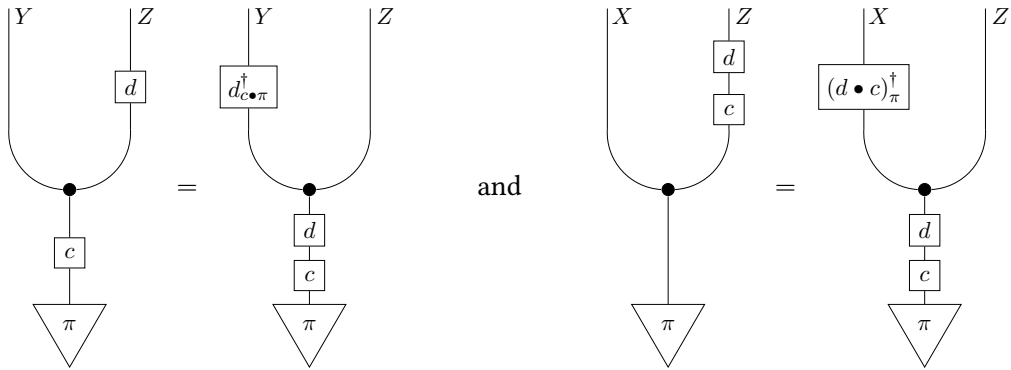
    \[
    \scalebox{0.8}{\begin{tikzpicture}
	\begin{pgfonlayer}{nodelayer}
		\node [style=black copier] (1) at (0, 0) {};
		\node [style=none] (2) at (-2, 6) {};
		\node [style=arrow box] (3) at (2, 3.5) {$d$};
		\node [style=none] (4) at (2, 6) {};
		\node [style=none] (5) at (-2, 2) {};
		\node [style=state180] (6) at (0, -4.5) {$\pi$};
		\node [style=arrow box] (8) at (0, -2) {$c$};
		\node [style=none] (9) at (-1.5, 5.75) {$Y$};
		\node [style=none] (10) at (2.5, 5.75) {$Z$};
		\node [style=none] (11) at (2, 2) {};
	\end{pgfonlayer}
	\begin{pgfonlayer}{edgelayer}
		\draw [bend left=45] (1) to (5.center);
		\draw (5.center) to (2.center);
		\draw (8) to (1);
		\draw (3) to (4.center);
		\draw (6) to (8);
		\draw [bend right=45] (1) to (11.center);
		\draw (11.center) to (3);
	\end{pgfonlayer}
\end{tikzpicture}}
    \ = \ 
    \scalebox{0.8}{\begin{tikzpicture}
	\begin{pgfonlayer}{nodelayer}
		\node [style=black copier] (0) at (0, 0) {};
		\node [style=none] (1) at (2, 6) {};
		\node [style=arrow box] (2) at (-2, 3.5) {$d^\dag_{c \bullet \pi}$};
		\node [style=none] (3) at (-2, 6) {};
		\node [style=none] (4) at (2, 2) {};
		\node [style=state180] (5) at (0, -4.5) {$\pi$};
		\node [style=arrow box] (6) at (0, -2.75) {$c$};
		\node [style=arrow box] (7) at (0, -1.25) {$d$};
		\node [style=none] (8) at (-1.5, 5.75) {$Y$};
		\node [style=none] (9) at (2.5, 5.75) {$Z$};
		\node [style=none] (10) at (-2, 2) {};
	\end{pgfonlayer}
	\begin{pgfonlayer}{edgelayer}
		\draw [bend right=45] (0) to (4.center);
		\draw (4.center) to (1.center);
		\draw (2) to (3.center);
		\draw (5) to (6);
		\draw (6) to (7);
		\draw (7) to (0);
		\draw [bend left=45] (0) to (10.center);
		\draw (10.center) to (2);
	\end{pgfonlayer}
\end{tikzpicture}}
    \hspace{0.06\linewidth}
    \text{and}
    \hspace{0.06\linewidth}
    \scalebox{0.8}{\begin{tikzpicture}
	\begin{pgfonlayer}{nodelayer}
		\node [style=black copier] (1) at (0, 0) {};
		\node [style=none] (2) at (-2, 6) {};
		\node [style=arrow box] (3) at (2, 2.75) {$c$};
		\node [style=none] (4) at (2, 6) {};
		\node [style=none] (5) at (-2, 2) {};
		\node [style=state180] (6) at (0, -4.5) {$\pi$};
		\node [style=arrow box] (7) at (2, 4.5) {$d$};
		\node [style=none] (8) at (-1.5, 5.75) {$X$};
		\node [style=none] (9) at (2.5, 5.75) {$Z$};
		\node [style=none] (10) at (2, 2) {};
	\end{pgfonlayer}
	\begin{pgfonlayer}{edgelayer}
		\draw [bend left=45] (1) to (5.center);
		\draw (5.center) to (2.center);
		\draw (6) to (1);
		\draw (3) to (7);
		\draw (7) to (4.center);
		\draw [bend right=45] (1) to (10.center);
		\draw (10.center) to (3);
	\end{pgfonlayer}
\end{tikzpicture}}
    \ = \ 
    \scalebox{0.8}{\begin{tikzpicture}
	\begin{pgfonlayer}{nodelayer}
		\node [style=black copier] (0) at (0, 0) {};
		\node [style=none] (1) at (2, 6) {};
		\node [style=arrow box] (2) at (-2, 3.5) {${(d \bullet c)}^\dag_\pi$};
		\node [style=none] (3) at (-2, 6) {};
		\node [style=none] (4) at (2, 2) {};
		\node [style=state180] (5) at (0, -4.5) {$\pi$};
		\node [style=arrow box] (6) at (0, -2.75) {$c$};
		\node [style=arrow box] (7) at (0, -1.25) {$d$};
		\node [style=none] (8) at (-1.5, 5.75) {$X$};
		\node [style=none] (9) at (2.5, 5.75) {$Z$};
		\node [style=none] (10) at (-2, 2) {};
	\end{pgfonlayer}
	\begin{pgfonlayer}{edgelayer}
		\draw [bend right=45] (0) to (4.center);
		\draw (4.center) to (1.center);
		\draw (2) to (3.center);
		\draw (5) to (6);
		\draw (6) to (7);
		\draw (7) to (0);
		\draw [bend left=45] (0) to (10.center);
		\draw (10.center) to (2);
	\end{pgfonlayer}
\end{tikzpicture}}
    \]

    The lens composite of these Bayesian inverses has the form \(c^\dag_\pi
    \klcirc d^\dag_{c \klcirc \pi} : Z \klto X\), so to establish the result it
    suffices to show that
    \[
    \scalebox{0.8}
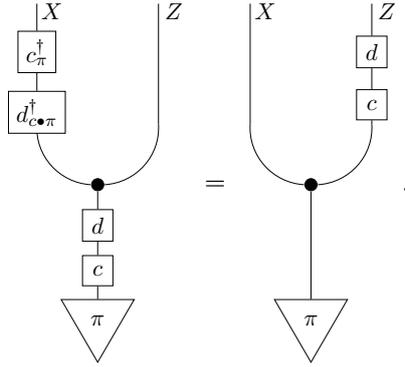
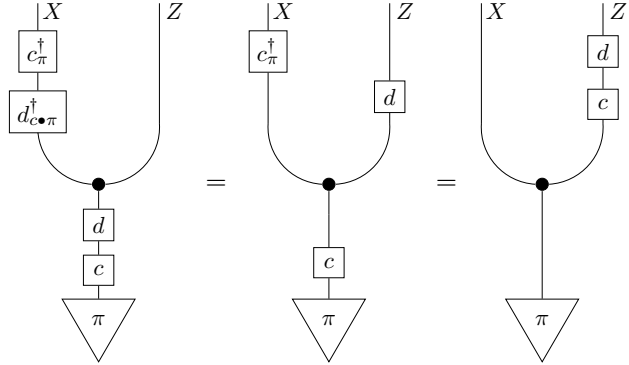
{\begin{tikzpicture}
	\begin{pgfonlayer}{nodelayer}
		\node [style=black copier] (0) at (0, 0) {};
		\node [style=none] (1) at (2, 6) {};
		\node [style=arrow box] (2) at (-2, 2.5) {$d^\dag_{c \bullet \pi}$};
		\node [style=none] (3) at (-2, 6) {};
		\node [style=none] (4) at (2, 2) {};
		\node [style=state180] (5) at (0, -4.5) {$\pi$};
		\node [style=arrow box] (6) at (0, -2.75) {$c$};
		\node [style=arrow box] (7) at (0, -1.25) {$d$};
		\node [style=arrow box] (8) at (-2, 4.5) {$c^\dag_\pi$};
		\node [style=none] (9) at (-1.5, 5.75) {$X$};
		\node [style=none] (10) at (2.5, 5.75) {$Z$};
		\node [style=none] (11) at (-2, 2) {};
	\end{pgfonlayer}
	\begin{pgfonlayer}{edgelayer}
		\draw [bend right=45] (0) to (4.center);
		\draw (4.center) to (1.center);
		\draw (5) to (6);
		\draw (6) to (7);
		\draw (7) to (0);
		\draw (2) to (8);
		\draw (8) to (3.center);
		\draw [bend left=45] (0) to (11.center);
		\draw (11.center) to (2);
	\end{pgfonlayer}
\end{tikzpicture}}
    \ = \ 
    \scalebox{0.8}{\begin{tikzpicture}
	\begin{pgfonlayer}{nodelayer}
		\node [style=black copier] (1) at (0, 0) {};
		\node [style=none] (2) at (-2, 6) {};
		\node [style=arrow box] (3) at (2, 2.75) {$c$};
		\node [style=none] (4) at (2, 6) {};
		\node [style=none] (5) at (-2, 2) {};
		\node [style=state180] (6) at (0, -4.5) {$\pi$};
		\node [style=arrow box] (7) at (2, 4.5) {$d$};
		\node [style=none] (8) at (-1.5, 5.75) {$X$};
		\node [style=none] (9) at (2.5, 5.75) {$Z$};
		\node [style=none] (10) at (2, 2) {};
	\end{pgfonlayer}
	\begin{pgfonlayer}{edgelayer}
		\draw [bend left=45] (1) to (5.center);
		\draw (5.center) to (2.center);
		\draw (6) to (1);
		\draw (3) to (7);
		\draw (7) to (4.center);
		\draw [bend right=45] (1) to (10.center);
		\draw (10.center) to (3);
	\end{pgfonlayer}
\end{tikzpicture}} .
    \]

    We have
    \[
    \scalebox{0.8}{\begin{tikzpicture}
	\begin{pgfonlayer}{nodelayer}
		\node [style=black copier] (0) at (0, 0) {};
		\node [style=none] (1) at (2, 6) {};
		\node [style=arrow box] (2) at (-2, 2.5) {$d^\dag_{c \bullet \pi}$};
		\node [style=none] (3) at (-2, 6) {};
		\node [style=none] (4) at (2, 2) {};
		\node [style=state180] (5) at (0, -4.5) {$\pi$};
		\node [style=arrow box] (6) at (0, -2.75) {$c$};
		\node [style=arrow box] (7) at (0, -1.25) {$d$};
		\node [style=arrow box] (8) at (-2, 4.5) {$c^\dag_\pi$};
		\node [style=none] (9) at (-1.5, 5.75) {$X$};
		\node [style=none] (10) at (2.5, 5.75) {$Z$};
		\node [style=none] (11) at (-2, 2) {};
	\end{pgfonlayer}
	\begin{pgfonlayer}{edgelayer}
		\draw [bend right=45] (0) to (4.center);
		\draw (4.center) to (1.center);
		\draw (5) to (6);
		\draw (6) to (7);
		\draw (7) to (0);
		\draw (2) to (8);
		\draw (8) to (3.center);
		\draw [bend left=45] (0) to (11.center);
		\draw (11.center) to (2);
	\end{pgfonlayer}
\end{tikzpicture}}
    \ = \ 
    \scalebox{0.8}{\begin{tikzpicture}
	\begin{pgfonlayer}{nodelayer}
		\node [style=black copier] (0) at (0, 0) {};
		\node [style=none] (1) at (2, 6) {};
		\node [style=none] (2) at (-2, 2) {};
		\node [style=none] (3) at (-2, 6) {};
		\node [style=arrow box] (4) at (2, 3) {$d$};
		\node [style=state180] (5) at (0, -4.5) {$\pi$};
		\node [style=arrow box] (6) at (0, -2.5) {$c$};
		\node [style=none] (7) at (0, -1) {};
		\node [style=arrow box] (8) at (-2, 4.5) {$c^\dag_\pi$};
		\node [style=none] (9) at (-1.5, 5.75) {$X$};
		\node [style=none] (10) at (2.5, 5.75) {$Z$};
		\node [style=none] (11) at (2, 2) {};
	\end{pgfonlayer}
	\begin{pgfonlayer}{edgelayer}
		\draw (4) to (1.center);
		\draw [bend left=45] (0) to (2.center);
		\draw (5) to (6);
		\draw (6) to (7.center);
		\draw (7.center) to (0);
		\draw (2.center) to (8);
		\draw (8) to (3.center);
		\draw [bend right=45] (0) to (11.center);
		\draw (11.center) to (4);
	\end{pgfonlayer}
\end{tikzpicture}}
    \ = \ 
    \scalebox{0.8}{\begin{tikzpicture}
	\begin{pgfonlayer}{nodelayer}
		\node [style=black copier] (1) at (0, 0) {};
		\node [style=none] (2) at (-2, 6) {};
		\node [style=arrow box] (3) at (2, 2.75) {$c$};
		\node [style=none] (4) at (2, 6) {};
		\node [style=none] (5) at (-2, 2) {};
		\node [style=state180] (6) at (0, -4.5) {$\pi$};
		\node [style=arrow box] (7) at (2, 4.5) {$d$};
		\node [style=none] (8) at (-1.5, 5.75) {$X$};
		\node [style=none] (9) at (2.5, 5.75) {$Z$};
		\node [style=none] (10) at (2, 2) {};
	\end{pgfonlayer}
	\begin{pgfonlayer}{edgelayer}
		\draw [bend left=45] (1) to (5.center);
		\draw (5.center) to (2.center);
		\draw (6) to (1);
		\draw (3) to (7);
		\draw (7) to (4.center);
		\draw [bend right=45] (1) to (10.center);
		\draw (10.center) to (3);
	\end{pgfonlayer}
\end{tikzpicture}}
    \]
    where the first obtains because \(d^\dag_{c \klcirc \pi}\) is the Bayesian
    inverse of \(d\) with respect to \(c \klcirc \pi\), and the second because
    \(c^\dag_\pi\) is the Bayesian inverse of \(c\) with respect to
    \(\pi\). Hence, \(c^\dag_\pi \klcirc d^\dag_{c \klcirc \pi}\) and \((d
    \klcirc c)^\dag_\pi\) are both Bayesian inversions of \(d \klcirc c\) with
    respect to \(\pi\). Since Bayesian inversions are almost-equal (Prop.
    \ref{prop:bayes-almost-equal}), we have \(c^\dag_\pi \klcirc d^\dag_{c
      \klcirc \pi} \overset{d \klcirc c \klcirc \pi}{\sim} (d \klcirc
    c)^\dag_\pi\), as required.
\end{proof}
\end{thm}

\begin{rmk}
  Note that, in the context of finitely-supported probability (\textit{e.g.}, in
  \(\Kl(\Da)\), where \(\Da\) is the finitely-supported probability distribution
  monad), almost-equality coincides with simple equality, and so Bayesian
  inversions are then just equal.
\end{rmk}

\begin{rmk} \label{rmk:lens-laws}
  Lenses were originally studied in the context of database systems
  \parencite{Bohannon2006Relational}, where one thinks of the forward channel as
  `viewing' a record in a database, and the backward channel as `updating' a
  record by taking a record and a new piece of data and returning the updated
  record. In this context, lenses have often been subject to additional axioms
  characterizing well-behavedness; for example, that updating a record with some
  data is idempotent (the `put-put' law). Bayesian lenses do not in general
  satisfy these laws, and nor even do exact Bayesian lenses. This is because
  Bayesian updating mixes information in the prior state (the `record') with the
  observation (the `data'), rather than replacing the prior information
  outright. We refer the reader to our preprint
  \parencite[{\S}6]{Smithe2020Bayesian} for a more detailed discussion of this
  situation.
\end{rmk}

\subsection{Parameterized Bayesian lenses} \label{sec:buco-para-blens}

Bayesian lenses for which the component channels are equipped with parameters
will play an important role in certain applications: an example which we will
meet in the next section is the \textit{variational autoencoder}, a neural
network architecture originally developed for machine learning and which is well
described by a particular class of parameterized statistical games. Since
\(\BLens{C}\) is a monoidal category, it induces a self-parameterization
\(\para(\BLens{C})\) by Proposition \ref{prop:para-self}, which is sufficient
for the purposes of this paper. Therefore, in this section, we summarize the
resulting structure for later reference.

\paragraph{0-cells and 1-cells}

The 0-cells of the bicategory \(\para(\BLens{C})\) are pairs \((X,A)\) of
objects in \(\cat{C}\). The 1-cells \((c,c^\dag) : (X,A) \xlensto{(\Omega,
  \Theta)} (Y,B)\) are Bayesian lenses \((c,c^\dag) : (\Omega \otimes X, \Theta
\otimes A) \lensto (Y, B)\). The forwards component \(c\) is a channel \(\Omega
\otimes X \klto Y\) in \(\cat{C}\), and hence also a parameterized channel \(X
\xklto{\Omega} Y\) in \(\para(\cat{C})\); we often think of this channel as
representing a system's model of the process by which observations (of
type~\(Y\)) are generated from causes (of type \(X\)), with the parameters (of
type \(\Omega\)) representing the system's beliefs about the structure of this
generative process.

Conversely, the backwards component \(c^\dag\) is a state-dependent channel \(B
\xklto{\Omega \otimes X} \Theta \otimes A\), which means a \(\Cat{V}\)-morphism
\(\cat{C}(I, \Omega \otimes X) \to \cat{C}(B, \Theta \otimes A)\). This is a
generalized Bayesian inversion which takes a state (or `prior' belief) jointly
over parameters and causes \(\Theta \otimes X\) and returns a channel \(B \klto
\Theta \otimes A\) that we think of as taking an observation (of possibly
different type \(B\)) and returning an updated joint belief about parameters and
causes (of possibly different types \(\Theta\) and \(A\)). Note that this
`update' channel is not itself parameterized: the parameterization of the
inversion is mediated through \(\Omega\), with \(\Theta\) being the type of
\textit{updated} parameters.

\paragraph{Sequential composition}

Given parameterized Bayesian lenses \((c,c^\dag) : (X,A)
\xlensto{(\Omega,\Theta)} (Y,B)\) and \((d,d^\dag) : (Y,B)
\xlensto{(\Omega',\Theta')} (Z,C)\), their composite \((d,d^\dag)\lenscirc
(c,c^\dag)\) is defined as the following morphism in \(\BLens{C}\):
\[
(\Omega' \otimes \Omega \otimes X, \Theta' \otimes
\Theta \otimes A) \xlensto{(\id, \id) \otimes (c,c^\dag)} (\Omega' \otimes Y,
\Theta' \otimes B) \xlensto{(d,d^\dag)} \, (Z, C) .
\]

\paragraph{Parallel composition}

Given parameterized Bayesian lenses \((c,c^\dag) : (X,A)
\xlensto{(\Omega,\Theta)} (Y,B)\) and \((d,d^\dag) : (X',A')
\xlensto{(\Omega',\Theta')} (Y',B')\), their tensor \((c,c^\dag) \otimes
(d,d^\dag)\) is defined as the following morphism in \(\BLens{C}\):
\[
\vlens{\Omega \otimes \Omega' \otimes X \otimes X'}{\Theta \otimes \Theta' \otimes A \otimes A'}
\xlensto{\sim}
\vlens{\Omega \otimes X \otimes \Omega' \otimes X'}{\Theta \otimes A \otimes \Theta' \otimes A'}
\xlensto{\vlens*{c}{\,c^\dag} \otimes \vlens*{d}{\,d^\dag}}
\vlens{Y \otimes Y'}{B \otimes B'}
\]
where we have written the pairs vertically, so that \(\vlens*{X}{A} := (X,A)\).

\paragraph{Reparameterization}

Given 1-cells \((f,f^\dag) : (X,A) \xlensto{(\Omega, \Theta)} (Y,B)\) and
\((g,g^\dag) : (X,A) \xlensto{(\Omega', \Theta')} (Y,B)\), the 2-cells \(\alpha
: (f,f^\dag) \Rightarrow (g,g^\dag)\) of \(\para(\BLens{C})\) are Bayesian
lenses \(\alpha : (\Omega,\Theta) \lensto (\Omega',\Theta')\) such that
\((f,f^\dag) = (g,g^\dag) \lenscirc (\alpha \otimes \id_{(X,A)})\). We can think
of the 2-cells as \textit{reparameterizations}, or higher-order processes that
predict the parameters on the basis of yet more abstract data.

\paragraph{Higher structure}

The \(\para\) construction turns a (monoidal) category into a (monoidal)
bicategory, thereby ``adding a dimension''. We can consider morphisms in a
parameterized category (as well as morphisms in \(\Fun{Stat}\); see Remark
\ref{rmk:stat-para-prox}) as processes by which processes are chosen, and
reparameterizations witness the factorization of these choice processes. In many
cybernetic and statistical situations, it is of interest to consider adding more
than just one extra dimension: that is, we may be interested in the processes by
which these choice processes are chosen, and, as in the Bayesian setting, we may
be interested in improving the performance of these `meta-processes'. More
concretely, in statistics, we may wish to describe meta-learning algorithms
(such as ``learning to learn''), or in neuroscience, we may wish to describe how
neuromodulation affects synaptic plasticity (which in turn affects the
generation of action potentials).

Since \(\para(\BLens{C})\) is itself monoidal, one has a further
parameterization \(\para(\para(\BLens{C}))\), and indeed \(\para\) has a monad
structure whose multiplication collapses a doubly-parameterized morphism to a
singly-parameterized one \parencite[Prop. 3]{Capucci2021Towards}. However, the
full structure that emerges when iterating this construction \textit{ad
  infinitum} is not yet well understood\footnote{At present, we believe the
resulting structure may have an opetopic shape, and constitute something like an
`\(\infty\)-fibration'.}. A similar structure appears when one considers
arbitrarily `nested' dynamical systems (such as cells in an organism, organisms
in societies, and societies in ecosystems; or the ownership structure of a
complex modern economy); see \textcite[Remark 2.3]{Smithe2021Some}. Such
cybernetic systems will be treated in a future chapter of the present series.

\section{Statistical games} \label{sec:sgame}

The Bayesian lenses of Theorem \ref{thm:buco} are exact, but most physically
realistic cybernetic systems do not compute exact inversions: the inversion of a
channel $c$ with respect to a prior $\pi$ generally involves evaluating the
composite $c \klcirc \pi$, which is typically computationally costly;
consequently, realistic systems typically instantiate approximate Bayesian
lenses. Fortunately, as a consequence of Theorem \ref{thm:buco}, an approximate
inversion of a composite channel will be approximately equal to this lens
composite; conversely, the lens composite of approximate inversions will be
approximately equal to the exact inversion of the corresponding composite.

We can thus approximate the inversion of a composite channel by the lens
composite of approximations to the inversions of the components. But what do we
mean by ``approximate''? There is often substantial freedom in the choice of
approximation scheme---often manifest as some form of parameterization---and in
typical situations, the `fitness' of a particular scheme will be
context-dependent. We think of Bayesian lenses as representing `open'
statistical systems, in interaction with some environment or `context', and so
the fitness of a particular lens then depends not only on the configuration of
the system (\emph{i.e.}, the choice of lens), but also on the suitability of its
configuration for the environment.

In this section, we introduce statistical games in order to quantify this
context-dependent fitness: a statistical game will be a Bayesian lens paired
with a contextual fitness function. The fitness function measures how well the
lens performs in each context, and the ``aim of the game'' is then to choose the
lens or context (depending on your perspective) that somehow optimizes the
fitness. In order to define statistical games, we first therefore define our
notion of context, following compositional game theory
\parencite{Ghani2016Compositional,Bolt2019Bayesian}: a context is ``everything
required to make an open system closed''. In the next section, we exemplify
statistical games by formalizing a number of classic and not-so-classic problems
in statistics.

\subsection{Contexts for Bayesian lenses}

Our first step is to define the notion of \textit{simple context}, with which we will be able to ``close off'' a lens with respect to sequential composition $\lenscirc$, and thereby define our categories of statistical games.

\begin{defn}
  A \textit{simple context} for a Bayesian lens over $\cat{C}$ is an element of the $\Cat{V}$-profunctor $\BLCtx{\cat{C}}$ defined by
  \[
  \begin{matrix*}[c]
    \BLens{\cat{C}} & \times & \BLens{\cat{C}}\op &   \to   & \Cat{V} \\
      {-}   & \times &     {= }   & \mapsto & \BLens{\cat{C}}((I,I),{-})\times\BLens{\cat{C}}({=},(I,I)) \, . \\
  \end{matrix*}
  \]
  If the lens is $(A,S)\lensto(B,T)$, its object of simple contexts is
  \[ \BLens{\cat{C}}\bigl((I,I),(A,S)\bigr)\times\BLens{\cat{C}}\bigl((B,T),(I,I)\bigr) \, . \]
  If we denote the lens by $f$, then we can denote this object of simple contexts by
  \[ \ctx(f) := \BLCtx{\cat{C}}\bigl((A,S),(B,T)\bigr) \, . \]
\end{defn}

When the monoidal unit $I$ is terminal in $\cat{C}$, or equivalently when $\cat{C}$ is semicartesian, the object of simple contexts acquires a simplified (and intuitive) form.

\begin{prop} \label{prop:ctx-terminal}
  When \(I\) is terminal in \(\cat{C}\), we have
  \[
  \BLCtx{C}\big((X,A),(Y,B)\big) \cong \cat{C}(I,X)\times\Cat{V}\bigl(\cat{C}(I,B),\cat{C}(I,Y)\bigr) \, .
  \]
  \begin{proof}
    A straightforward calculation which we omit: use causality (Definition
    \ref{def:causal}).
  \end{proof}
\end{prop}

\begin{rmk}
  Proposition \ref{prop:ctx-terminal} means that, when \(\cat{C}\) is a Markov category (such as \(\Kl(\Giry)\)), a context for a Bayesian lens consists of a `prior' on the domain of the forwards channel and a `continuation':
  a \(\Cat{V}\)-morphism (such as a function) which takes the output of the forwards channel (possibly a `prediction') and returns an observation for the update map.
  We think of the continuation as encoding the response of the environment given the prediction.
\end{rmk}

In order to define the sequential composition of statitistical games, we will need to construct, from the context for a composite lens, contexts for each factor of the composite. The functoriality of \(\BLCtx{C}\) guarantees the existence of such \textit{local contexts}.

\begin{defn}
  Given another lens $g : (B,T)\lensto(C,U)$ and a simple context for their composite $g\lenscirc f : (A,S)\lensto(C,Y)$, then we can obtain a \textit{1-local} context for $f$ by the action of the profunctor
  \[
  \BLCtx{\cat{C}}\bigl((A,S),g\bigr) : \BLCtx{\cat{C}}\bigl((A,S),(C,U)\bigr) \to \BLCtx{\cat{C}}\bigl((A,S),(B,T)\bigr)
  \]
  and we can obtain a simple context for $g$ similarly:
  \[
  \BLCtx{\cat{C}}\bigl(f,(C,U)\bigr) : \BLCtx{\cat{C}}\bigl((A,S),(C,U)\bigr) \to \BLCtx{\cat{C}}\bigl((B,T),(C,U)\bigr) \, .
  \]
  Note that we can write $g^* : \ctx(g\lenscirc f) \to \ctx(f)$ for the former action and $f_* : \ctx(g\lenscirc f) \to \ctx(g)$ for the latter, since $g^*$ acts by precomposition (pullback) and $f_*$ acts by postcomposition (pushforwards).
  Note also that, given $h\lenscirc g\lenscirc f$, we have $h^*\circ f_* = f_*\circ h^* : \ctx(h\lenscirc g\lenscirc f) \to \ctx(g)$ by the associativity of composition.
\end{defn}

\begin{rmk}
  The local contexts of the preceding definition are local with respect to sequential composition of lenses.
  If we view the monoidal category $\BLens{\cat{C}}$ as a one-object bicategory, then sequential composition is composition of 1-cells, which explains their formal naming as \textit{1-local} contexts.
\end{rmk}

To lift the monoidal structure of $\BLens{\cat{C}}$ to our categories of statistical games, we will similarly need `2-local' contexts:
from the one-object bicategory perspective, these are local contexts with respect to 2-cell composition.
By analogy with the 1-local case, this means exhibiting maps of the form $\ctx(f\otimes f')\to\ctx(f)$ and $\ctx(f\otimes f')\to\ctx(f')$ which give the `left' and `right' local contexts for a parallel pair of lenses.

Without requiring extra structure from $\cat{C}$ or $\BLens{C}$\footnote{
An alternative is to ask for $\BLens{C}$ to be equipped with a natural family of `discarding' morphisms $(A,S)\to(I,I)$.
This in turn means asking for canonical states $I\to S$, which in general we do not have:
if we are working with stochastic channels, we could obtain such canonical states by allowing the channels to emit subdistributions, and letting the canonical states be given by those which assign $0$ density everywhere.
But even though we can make the types check this way, the semantics are not quite what we want: rather, we seek to allow information to ``pass in parallel''.},
we first need to pass from simplex contexts to \textit{complex} ones: a 2-local context should allow for information to pass alongside a lens, through a process that is parallel to it.
Formally, we adjoin an object to the domain and codomain of the context, and quotient by the rule that, if there is any process that `fills the hole' represented by the adjoined objects, then we consider the objects equivalent:
this allows us to forget about the contextually parallel processes, and keep track only of the type of information that flows. Such adjoining-and-quotienting gives us the following coend\footnote{
For more information about \textit{(co)ends} and \textit{(co)end calculus}, we refer the reader to \textcite{Loregian2021Coend}.}
formula defining complex contexts.

\begin{defn}
  We define a \textit{complex context} to be an element of the profunctor
  \[
  \BLCCtx{C}\bigl((A,S), (B,T)\bigr) := \int^{(M,N):\BLens{C}} \BLCtx{C}\bigl((M,N)\otimes(A,S), (M,N)\otimes(B,T)\bigr) \, .
  \]
  If $f : (A,S)\lensto(B,T)$ is a lens, then we will write $\cctx(f) := \BLCCtx{C}\bigl((A,S),(B,T)\bigr)$.
  As in the case of simple contexts, we have `projection' maps $g^* : \cctx(g\lenscirc f)\to\cctx(f)$ and $f_* : \cctx(g\lenscirc f)\to\cctx(g)$; these are defined similarly.
  We will call the adjoined object, here denoted $(M,N)$, the \textit{residual}.
  There is a canonical inclusion $\ctx(f)\hookrightarrow\cctx(f)$ given by adjoining the trivial residual $(I,I)$ to each simple context.
\end{defn}

We immediately have the following corollary of Proposition \ref{prop:ctx-terminal}:

\begin{cor} \label{cor:comp-ctx-term}
  When $I$ is terminal in $\cat{C}$,
  \[
  \BLCCtx{C}\bigl((A,S),(B,T)\bigr) \cong \int^{(M,N)} \cat{C}(I, M\otimes A) \times \Cat{V}\bigl(\cat{C}(I, M\otimes B), \cat{C}(I, N\otimes T)\bigr) \, .
  \]
\end{cor}

\begin{rmk}
  Since the coend denotes a quotient, its elements are equivalence classes.
  The preceding corollary says that, when $I$ is terminal, an equivalence class of contexts is represented by a choice of residual $(M,N)$, a prior on $M\otimes A$ in $\cat{C}$, and a continuation $\cat{C}(I, M\otimes B)\to\cat{C}(I, N\otimes T)$ in $\Cat{V}$.

  Of course, when $I$ is not terminal, then the definition says that a general complex context for a lens $(A,S)\lensto(B,T)$ is an equivalence class represented by: as before,
  a choice of residual $(M,N)$,
  a prior on $M\otimes A$,
  and a continuation $\cat{C}(I,M\otimes B)\to\cat{C}(I,N\otimes T)$;
  as well as an `effect' $M\otimes B\klto I$ in $\cat{C}$,
  and what we might call a `vector' $\cat{C}(I,I)\to\cat{C}(N\otimes S, I)$ in $\Cat{V}$.
  The effect and vector measure the environment's `internal response' to the lens' outputs (as opposed to representing the environment's feedback to the lens).
\end{rmk}

Using complex contexts, it is easy to define the `projections' that give 2-local contexts.

\begin{defn} \label{def:2local-ctx}
  Given lenses $f:\Phi\lensto\Psi$ and $f':\Phi'\lensto\Psi'$ and a complex context for their tensor $f\otimes f'$, the \textit{left 2-local context} is the complex context for $f$ given by
  \begin{align*}
    \pi_f : \; & \cctx(f\otimes f') = \BLCCtx{C}\bigl(\Phi\otimes\Phi',\Psi\otimes\Psi'\bigr) \\
    & = \int^{\Theta:\BLens{C}} \BLCtx{C}\bigl(\Theta\otimes\Phi\otimes\Phi',\Theta\otimes\Psi\otimes\Psi'\bigr) \\
    & \xto{\sim} \int^{\Theta:\BLens{C}} \BLCtx{C}\bigl(\Theta\otimes\Phi'\otimes\Phi,\Theta\otimes\Psi'\otimes\Psi\bigr) \\
    & = \int^{\Theta:\BLens{C}} \BLens{C}\bigl((I,I),\Theta\otimes\Phi'\otimes\Phi\bigr)\times\BLens{C}\bigl(\Theta\otimes\Psi'\otimes\Psi,(I,I)\bigr) \\
    & \xto{f'} \int^{\Theta:\BLens{C}} \BLens{C}\bigl((I,I),\Theta\otimes\Psi'\otimes\Phi\bigr)\times\BLens{C}\bigl(\Theta\otimes\Psi'\otimes\Psi,(I,I)\bigr) \\
    & = \int^{\Theta:\BLens{C}} \BLCtx{C}\bigl(\Theta\otimes\Psi'\otimes\Phi,\Theta\otimes\Psi'\otimes\Psi\bigr) \\
    & \hookrightarrow \int^{\Theta':\BLens{C}} \BLCtx{C}\bigl(\Theta'\otimes\Phi,\Theta'\otimes\Psi\bigr) = \cctx(f) \, .
  \end{align*}
  The equalities here are just given by expanding and contracting definitions; the isomorphism uses the symmetry of $\otimes$ in $\cat{C}$; the arrow marked $f'$ is given by composing $\id_\Theta\otimes f'\otimes\id_\Phi$ after the `prior' part of the context; and the inclusion is given by collecting the tensor of $\Theta$ and $\Psi'$ together into the residual. These steps formalize the idea of filling the right-hand hole of the complex context with $f'$ to obtain a local context for $f$.

  The \textit{right 2-local context} is the complex context for $f'$ obtained similarly:
  \begin{align*}
    \pi_{f'} : \; & \cctx(f\otimes f') = \BLCCtx{C}\bigl(\Phi\otimes\Phi',\Psi\otimes\Psi'\bigr) \\
    & = \int^{\Theta:\BLens{C}} \BLCtx{C}\bigl(\Theta\otimes\Phi\otimes\Phi',\Theta\otimes\Psi\otimes\Psi'\bigr) \\
    & = \int^{\Theta:\BLens{C}} \BLens{C}\bigl((I,I),\Theta\otimes\Phi\otimes\Phi'\bigr)\times\BLens{C}\bigl(\Theta\otimes\Psi\otimes\Psi',(I,I)\bigr) \\
    & \xto{f} \int^{\Theta:\BLens{C}} \BLens{C}\bigl((I,I),\Theta\otimes\Psi\otimes\Phi'\bigr)\times\BLens{C}\bigl(\Theta\otimes\Psi\otimes\Psi',(I,I)\bigr) \\
    & = \int^{\Theta:\BLens{C}} \BLCtx{C}\bigl(\Theta\otimes\Psi\otimes\Phi',\Theta\otimes\Psi\otimes\Psi'\bigr) \\
    & \hookrightarrow \int^{\Theta':\BLens{C}} \BLCtx{C}\bigl(\Theta'\otimes\Phi',\Theta'\otimes\Psi'\bigr) = \cctx(f')
  \end{align*}
  Note here that we do not need to use the symmetry, as both the residual and the `hole' (filled with $f$) are on the left of the tensor. (Strictly, we do not require symmetry of $\otimes$ for $\pi_f$ either, only a braiding.)
\end{defn}

\begin{rmk}
  It is possible to make the intuition of ``filling the left and right holes''
  more immediately precise, at the cost of introducing another language, by
  rendering the 2-local context functions in the graphical calculus of the
  monoidal bicategory of \(\Cat{V}\)-profunctors. We demonstrate how this works,
  making the hole-filling explicit, in Appendix \ref{sec:apdx-2ctx}.
\end{rmk}

Henceforth, when we say `context', we will mean `complex context'.

\subsection{Monoidal categories of statistical games}

We are now in a position to define monoidal categories of statistical games over \(\cat{C}\).
In typical examples, the fitness functions will be valued in the real numbers, but this is not necessary for the categorical definition;
instead, we allow the fitness functions to take values in an arbitrary monoid.

\begin{prop} \label{prop:cat-sgame}
  Let \(\cat{C}\) be a \(\Cat{V}\)-category admitting Bayesian inversion and let
  \((R, +, 0)\) be a monoid in \(\Cat{V}\). Then there is a category
  \(\SGame[R]{\cat{C}}\) whose objects are the objects of \(\BLens{C}\) and
  whose morphisms \((X, A) \to (Y, B)\) are \textit{statistical games}: pairs
  \((f, \phi)\) of a lens \(f : \BLens{C}\big((X, A), (Y, B)\big)\) and a
  \textit{fitness function} \(\phi : \cctx(f) \to R\).
  When \(R\) is the monoid of reals \(\rr\), then we just denote the category by
  \(\SGame{\cat{C}}\).
  \begin{proof}
    Suppose given statistical games \((f,\phi) : (X,A) \to (Y,B)\) and \((g,\psi) : (Y,B) \to (Z,C)\).
    We seek a composite game \((g,\psi) \circ (f,\phi) := (gf,\psi\phi) : (X,A) \to (Z,C)\).
    We have \(gf = g \lenscirc f\) by lens composition.
    The composite fitness function \(\psi\phi\) is obtained using the 1-local contexts by
    \[
    \psi\phi := \cctx(g\lenscirc f) \xto{\bcopier} \cctx(g\lenscirc f)\times\cctx(g\lenscirc f) \xto{(g^*, f_*)} \cctx(f)\times\cctx(g) \xto{(\phi, \psi)} R\times R \xto{+} R
    \]
    The identity game \((X,A) \to (X,A)\) is given by \((\id, 0)\), the pairing of the identity lens on \((X, A)\) with the unit \(0\) of the monoid \(R\).
    Associativity and unitality follow from those properties of lens composition, the coassociativity of copying $\bcopier$ in $\Cat{V}$, and the monoid laws of $R$.
  \end{proof}
\end{prop}

\begin{defn} \label{def:simp-sgame}
  We will write \(\SimpSGame{\cat{C}} \hookrightarrow \SGame{\cat{C}}\) for the
  full subcategory of \(\SGame{\cat{C}}\) defined on simple Bayesian lenses
  \((X,X) \lensto (Y, Y)\). As the case of simple lenses (Definition
  \ref{def:simp-lens}), we will eschew redundancy by writing the objects
  \((X,X)\) simply as \(X\).
\end{defn}

\begin{prop} \label{prop:sgame-smc}
  \(\SGame[R]{\cat{C}}\) inherits a monoidal structure \(\big(\otimes,(I,I)\big)\) from \(\big(\BLens{C},\otimes,(I,I)\big)\).
  When \((\cat{C}, \otimes, I)\) is furthermore symmetric monoidal and \((R, +, 0)\) is a commutative monoid, then the monoidal structure on \(\SGame[R]{\cat{C}}\) is
  symmetric.
  \begin{proof}
    The structure on objects and on the lens components of games is defined as in Proposition \ref{prop:blens-monoidal}.
    On fitness functions, we use the monoidal structure of \(R\):
    given games \((f,\phi) : (X,A) \to (Y,B)\) and \((f',\phi') : (X',A') \to (Y',B')\), we define the fitness function \(\phi \otimes \phi'\) of their tensor \((f,\phi) \otimes (f',\phi')\) as the composite
    \[
    \phi\otimes\phi' := \cctx(f\otimes f') \xto{\bcopier} \cctx(f\otimes f')\times\cctx(f\otimes f') \xto{(\pi_f,\pi_{f'})} \cctx(f)\times\cctx(f') \xto{(\phi,\phi')} R\times R \xto{+} R
    \]
    That is, we form the left and right 2-local contexts, compute the local fitnesses, and compose them using the monoidal operation in \(R\).
    Unitality and associativity follow from that of \(\otimes\) in \(\BLens{C}\) and \(+\) in \(R\).
  \end{proof}
\end{prop}

\begin{cor}
  Since each category \(\SGame[R]{\cat{C}}\) is thus monoidal, we obtain
  categories \(\para(\SGame[R]{\cat{C}})\) of \textit{parameterized statistical
    games} by Proposition \ref{prop:para-self}, which are themselves
  monoidal. The structure is as described for Bayesian lenses in Section
  \ref{sec:buco-para-blens}, with some minor additions to incorporate the
  fitness functions, the details of which we leave to the reader.
\end{cor}

\begin{rmk}[Parameters as strategies]
  A parameterized statistical game of type \((X,A) \xto{(\Omega,\Theta)} (Y,B)\)
  in \(\para(\SGame[R]{\cat{C}})\) is a statistical game \(\vlens*{\Omega
    \otimes X}{\Theta \otimes A} \to \vlens*{Y}{B}\) in \(\SGame[R]{\cat{C}}\);
  that is a pair of a Bayesian lens \(\vlens*{\Omega \otimes X}{\Theta \otimes
    A} \lensto \vlens*{Y}{B}\) and a fitness function
  \(\BLCCtx{C}\left(\vlens*{\Omega \otimes X}{\Theta \otimes A},
  \vlens*{Y}{B}\right) \to R\) in \(\Cat{V}\). If we fix a choice of parameter
  \(\omega : \Theta\) and discard the updated parameters in \(\Theta\)---that
  is, if we reparameterize along the 2-cell induced by the lens \((\omega,
  \ground) : (I,I) \lensto (\Omega,\Theta)\)---then we obtain an unparameterized
  statistical game \((X,A) \to (Y,B)\). In this way, we can think of the
  parameters of a parameterized statistical game as the \textit{strategies} by
  which the game is to be played: each parameter \(\omega : \Omega\) picks out a
  Bayesian lens, whose forwards channel we think of as a model by which the
  system predicts its observations and whose backwards channel describes how the
  system updates its beliefs. And, if we don't just discard them, then these
  updated beliefs may include updated parameters (of a possibly different type
  \(\Theta\)). A successful strategy (a good choice of parameter) for a
  statistical game is then one which optimizes the fitness function in the
  contexts of relevance to the system: we can think of these ``relevant
  contexts'' as something like the system's ecological niche.
\end{rmk}

\begin{rmk}[Multi-player games]
  In a later instalment of this series of papers, we will see how to compose the
  statistical games of multiple interacting agents, so that the game-playing
  metaphor becomes more visceral: the observations predicted by each system will
  then be generated by other systems, with each playing a game of optimal
  prediction. In this paper, however, we usually think of each statistical game
  as representing a single system's model of its environment (its context), even
  where the games at hand are themselves sequentially or parallelly
  composite. That is to say, our games here are fundamentally `two-player', with
  the two players being the system and the context.
\end{rmk}

\begin{rmk} \label{rmk:fit-glue}
  Both in the case of sequential of parallel composition of statistical games, the local fitnesses are computed independently and then summed.
  If a fitnesss function depends somehow on the residual, this might lead to `double-counting' the fitness of any overlapping factors of the residual.
  For our purposes, this assumption of `independent fitness' will suffice, and so we leave the question of gluing together correlated fitness functions for future work.
\end{rmk}

\section{Examples} \label{sec:examples}

In this section, we describe how a number of common concepts in statistics and
particularly statistical inference fit into the framework of statistical
games. We begin with the simple example of maximum likelihood estimation and
progressively generalize to include `variational'
\parencite{Blei2016Variational} methods such as the variational autoencoder
\parencite{Kingma2013Auto} and generalized variational inference
\parencite{Knoblauch2019Generalized}. Along the way, we introduce the concepts
of \textit{free energy} and \textit{evidence upper bound}.  We do not here
consider the algorithms by which statistical games may be played or optimized;
that is a matter for a subsequent paper in this series. Instead, we see
statistical games as providing an `algebra' for the compositional construction
of inference problems.

\begin{rmk}[The role of fitness functions]
  Before we introduce our first example, we note that the games here are classified by their fitness functions, with the choice of lens being somewhat incidental to the classification\footnote{
  This incidentality is lessened when we consider examples of \textit{parameterized} games, but even here the parameterization only induces something of a `sub'-classification; the main
  classification remains due to the fitness functions.}.
  We note furthermore that our fitness functions will tend to be of the form \(\E_{k\klcirc c\klcirc\pi}[f]\), where \((\pi, k)\) is a context for a lens, \(c\) is a channel, and \(f\) is an appropriately typed effect\footnote{
  The resulting `optimization-centric' perspective is in line with the aesthetic preference of \parencite{Knoblauch2019Generalized}, though we do not yet know what this alignment might signify; we are interested to find examples of a different flavour.}.
  This form hints at the existence of a compositional treatment of fitness functions, which seem roughly to be something like ``lens functionals''.
  We leave such a treatment, and its connection to Remark \ref{rmk:fit-glue}, to future work.
\end{rmk}

We first study the classic problem of maximum likelihood estimation, beginning by establishing an auxiliary results about contexts.

\begin{prop} \label{prop:mle-ctx}
  Let $\mathrm{I}$ denote the monoidal unit $(I,I)$ in $\BLens{C}$, and let $l : \mathrm{I}\lensto\Psi$ be a lens. Then
  \begin{align*}
    \cctx(l) & = \int^{\Theta:\BLens{C}} \BLens{C}(\mathrm{I}, \Theta\otimes\mathrm{I}) \times \BLens{C}(\Theta\otimes\Psi,\mathrm{I}) \\
    & \cong \int^{\Theta:\BLens{C}} \BLens{C}(\mathrm{I},\Theta)\times\BLens{C}(\Theta\otimes\Psi,\mathrm{I}) \\
    & \cong \BLens{C}(\mathrm{I}\otimes\Psi,\mathrm{I}) \\
    & \cong \BLens{C}(\Psi,\mathrm{I}) \, .
  \end{align*}
  Suppose $\Psi = (A,S)$. Then $\cctx(l) = \cat{C}(A,I)\times\Cat{V}\bigl(\cat{C}(I,A),\cat{C}(I,S)\bigr)$ by the definition of $\BLens{C}$.
  \begin{proof}
    The first and third isomorphisms hold by unitality of $\otimes$;
    the second holds by the Yoneda lemma (see \textcite[Prop. 2.2.1]{Loregian2021Coend} for the argument).
  \end{proof}
\end{prop}

\begin{ex}[Maximum likelihood] \label{ex:ml-game}
  When \(I\) is terminal in \(\cat{C}\), a Bayesian lens of the form \((I, I) \lensto (X, X)\) is determined by its forwards channel, which is simply a state \(\pi : I \klto X\).
  Following Proposition \ref{prop:mle-ctx}, and using that $I$ is terminal, a context for such a lens is given simply by a continuation $k : \cat{C}(I,X)\to\cat{C}(I,X)$ taking states on $X$ to states on $X$.
  A \textit{maximum likelihood game} is then any statistical game $\pi$ of the type \((I, I) \to (X, X)\) with fitness function \(\phi : \cctx(\pi)\to\rr\) given by
  \(\phi(k) = \E_{k(\pi)} \left[ p_\pi \right]\),
  where \(p_\pi\) is a density function for \(\pi\).
  More generally, we might consider \textit{maximum \(f\)-likelihood games} for monotone functions \(f : \rr \to \rr\), in which the fitness function is given by
  \(\phi(k) = \E_{k(\pi)} \left[ f \circ p_\pi \right]\).
  A typical choice here is \(f := \log\).
\end{ex}

In order that there may be some freedom to optimize the fitness function, one typically works in the parameterized category:
the aim of the game is then to choose the optimal parameter for the context, as quantified by the fitness function.
This gives us the notion of \textit{parameterized maximum likelihood} game;
but first, we define some simplifying notation.

\begin{notation}[Feedback]
  Let $I$ be terminal in $\cat{C}$ and consider a Bayesian lens $l = (l_1,l'):(A,S)\lensto(B,T)$, with a context represented by: a residual $(M,N)$; a prior $\pi:I\klto M\otimes A$; and a continuation $k:\cat{C}(I,M\otimes B)\to\cat{C}(I,N\otimes T)$.
  Write $\efb{\pi}{l}{k}$ to denote $k\bigl((\id_M\otimes\, l_1)\klcirc\pi\bigr)_T$ where $(\id_M\otimes\,l_1)\klcirc\pi$ is the map
  \[
  I\xklto{\pi}M\otimes A\xklto{\id_M\otimes\,l_1}M\otimes B
  \]
  and where $(-)_T$ denotes the projection (marginalization) onto $T$; here, by the channel $N\otimes T\klto T$.

  Note that $\efb{\pi}{l}{k}$ therefore has the type $I\klto T$ in $\cat{C}$: it encodes the environment's feedback in $T$ to the lens, given its output in $B$ and the context.
\end{notation}

\begin{ex}[Parameterized maximum likelihood] \label{ex:para-ml-game}
  A parameterized Bayesian lens \((I,I) \xlensto{(\Omega,\Theta)} (X,X)\) is equivalently a Bayesian lens \((\Omega,\Theta)\lensto(X,X)\), and hence given by a pair of a channel (or ``parameter-dependent state'') \(\Omega \klto X\) and a parameter-update \(X \xklto{\Omega} \Theta\).
  When $I$ is terminal in $\cat{C}$, a context for the lens is represented by $\pi:I\klto M\otimes\Omega$ and $k:\cat{C}(I,M\otimes X)\to\cat{C}(I,N\otimes X)$.
  We then define a \textit{parameterized maximum \(f\)-likelihood game} to be a parameterized statistical game of the form \(l = (l_1,l') : (I,I) \xlensto{(\Omega,\Theta)} (X,X)\) with fitness function \(\phi : \cctx(\pi) \to \rr\) given by
  \( \phi(\pi,k) = \E_{\efb{\pi}{l}{k}} \left[ f \circ p_{l\klcirc\pi_\Omega} \right] \) .

  Here, \(p_{l\klcirc\pi_\Omega}:X\to[0,\infty]\) is a density function for the composite channel \(l\klcirc\pi_\Omega\).
  In applications one often fixes a single choice of parameter $o$, with the marginal state \(\pi_\Omega\) then being a Dirac delta distribution on that choice.
  One then writes the density function $p_{l\klcirc\pi_\Omega}({-})$ as $p_l({-}|o)$ or $p_l\left({-}|\Omega=o\right)$.
\end{ex}

\begin{rmk} \label{rmk:ml-game}
  Recalling that we can think of probability density as a measure of the likelihood of an observation, we have the intuition that an ``optimal strategy'' (\textit{i.e.}, an optimal choice of lens or parameter) for a maximum likeihood game is one that maximizes the likelihood of the state obtained from the context, or in other words provides the ``best explanation'' of the data generated by the continuation.
\end{rmk}

Considering parameterized maximum likelihood games, which are equipped with
parameter-update maps, leads one to wonder how to optimize this `inferential'
backwards part of the game, and not just the `predictive' forwards part. Such
backwards optimization is approximate Bayesian inference.

\begin{ex}[Bayesian inference] \label{ex:bayes-game}
  Let \(D : \cat{C}(I,X) \times \cat{C}(I,X) \to \rr\) be a measure of
  divergence between states on \(X\). Then a (simple) \(D\)\textit{-Bayesian
    inference} game is a statistical game \((c,\phi) : (X,X) \to (Y,Y)\) with fitness
  function \(\phi : \cctx(c)\to\rr\) given by
  \(\phi(\pi,k) = \E_{y\sim\efb{\pi}{c}{k}} \left[ D\left(c'_\pi(y),
    c^\dag_\pi(y)\right) \right]\), where \(c = (c_1,c')\) constitutes the lens part
  of the game and \(c^\dag_\pi\) is the exact inversion of \(c_1\) with respect to
  \(\pi\).
\end{ex}

Note that we say that \(D\) is a ``measure of divergence between states on
\(X\)''. By this we mean any function of the given type with the semantical
interpretation that it acts like a distance measure between states. But this is
not to say that \(D\) is a metric or even pseudometric. One usually requires
that \(D(\pi, \pi') = 0 \iff \pi = \pi'\), but typical choices do not also
satisfy symmetry nor subadditivity. An important such typical choice is the
\textit{relative entropy} or \textit{Kullback-Leibler divergence}, denoted
\(D_{KL}\).

\begin{defn} \label{def:d-kl}
  The \textit{Kullback-Leibler divergence} \(D_{KL} : \cat{C}(I, X) \times
  \cat{C}(I, X) \to \rr\) is defined by
  \[
  D_{KL}(\alpha, \beta) := \E_{x \sim \alpha}[\log p (x)] - \E_{x \sim \alpha}[\log q (x)]
  \]
  where \(p\) and \(q\) are density functions corresponding to the states
  \(\alpha\) and \(\beta\).
\end{defn}

In many situations, computing the exact inversion \(c^\dag_\pi(x)\) is
costly, and so is computing the divergence \(D\left(c'_\pi(x),
c^\dag_\pi(x)\right)\). Consequently, approximate inversion schemes typically
either approximate the divergence (as in Monte Carlo methods), or they optimize
an upper bound on it (as in variational methods). In this section, we are
interested in different choices of fitness function, rather than the algorithms
by which the functions are exactly or approximately evaluated; hence we here
consider the latter `variational' choice, leaving the former for future work.

One widespread choice is to construct an upper bound on the divergence called
the \textit{free energy} or the \textit{evidence upper bound}.

\begin{defn}[\(D\)-free energy] \label{def:d-free-energy}
  Let \((\pi, c)\) be a generative model with \(c : X \klto Y\). Let \(p_c : Y
  \times X \to \rr_+\) and \(p_\pi : X \to \rr_+\) be density functions
  corresponding to \(c\) and \(\pi\). Let \(p_{c \klcirc \pi} : Y \to \rr_+\) be
  a density function for the composite \(c\klcirc\pi\).  Let \(c'_\pi\) be a
  channel \(Y \klto X\) that we take to be an approximation of the Bayesian
  inversion of \(c\) with respect to \(\pi\)and that admits a density function
  \(q : X \times Y \to \rr_+\). Finally, let \(D : \cat{C}(I, X) \times
  \cat{C}(I, X) \to \rr\) be a measure of divergence between states on \(X\).
  Then the \(D\)\textit{-free energy} of \(c'_\pi\) with respect to the
  generative model given an observation \(y:Y\) is the quantity
  \begin{equation} \label{eq:free-energy}
  \Fa_D(c'_\pi, c, \pi, y) := \E_{x \sim c'_\pi(y)} \left[ - \log p_c(y | x)  \right] + D\left(c'_\pi(y), \pi\right) \, .
  \end{equation}
  We will elide the dependence on the model when it is clear from the context,
  writing only \(\Fa_D(y)\).
\end{defn}

The \(D\)-free energy is an upper bound on \(D\) when \(D\) is the relative
entropy \(D_{KL}\), as we now show.

\begin{prop}[Evidence upper bound] \label{prop:eubo}
  The \(D_{KL}\)-free energy satisfies the following equality:
  \[
  \Fa_{D_{KL}}(y) = D_{KL}\left[ c'_\pi(y), c^\dag_\pi(y)\right] - \log p_{c \klcirc \pi}(y) = \E_{x \sim c'_\pi(y)} \left[ \log \frac{q(x|y)}{p_c(y|x) \cdot p_\pi(x)} \right]
  \]
  Since \(\log p_{c \klcirc \pi}(y)\) is always negative, the free energy is an
  upper bound on \(D_{KL}\left[ c'_\pi(y), c^\dag_\pi(y)\right]\), where
  \(c^\dag_\pi\) is the exact Bayesian inversion of the channel \(c\) with
  respect to the prior \(\pi\). Similarly, the free energy is an upper bound on
  the negative log-likelihood \(-\log p_{c \klcirc \pi}(y)\). Thinking of this
  latter quantity as a measure of the ``model evidence'' gives us the
  alternative name \textit{evidence upper bound} for the \(D_{KL}\)-free energy.
  \begin{proof}
    Let \(p_\omega : Y \times X \to \rr_+\) be the density function
    \(p_\omega(y, x) := p_c(y|x) \cdot p_\pi(x)\) corresponding to the joint
    distribution of the generative model \((\pi, c)\). We have the following
    equalities:
    \begin{align*}
      -\log p_{c \klcirc \pi}(y)
      &= \E_{x \sim c'_\pi(y)} \left[ - \log p_{c \klcirc \pi}(y) \right] \\
      &= \E_{x \sim c'_\pi(y)} \left[ - \log \frac{p_\omega(y,x)}{p_{c^\dag_\pi}(x|y)} \right]
      \qquad\text{(by Bayes' rule)} \\
      &= \E_{x \sim c'_\pi(y)} \left[ - \log \frac{p_\omega(y,x)}{q(x|y)} \frac{q(x|y)}{p_{c^\dag_\pi}(x|y)} \right] \\
      &= - \E_{x \sim c'_\pi(y)} \left[ \log \frac{p_\omega(y,x)}{q(x|y)} \right] - D_{KL}\left[ c'_\pi(y), c^\dag_\pi(y) \right]
    \end{align*}
  \end{proof}
\end{prop}

\begin{defn}
  We will call \(\Fa_{D_{KL}}\) the \textit{variational free energy}, or simply
  \textit{free energy}, and denote it by \(\Fa\) where this will not cause
  confusion. We will take the result of Proposition \ref{prop:eubo} as a
  definition of the variational free energy, writing
  \[
  \Fa(y) = \E_{x \sim c'_\pi(y)} \left[ \log \frac{q(x|y)}{p_c(y|x) \cdot p_\pi(x)} \right]
  \]
  where each term is defined as in Definitions \ref{def:d-free-energy} and \ref{def:d-kl}.
\end{defn}

\begin{rmk} \label{rmk:helmholtz}
  The name \textit{free energy} is due to an analogy with the Helmholtz free
  energy in thermodynamics, as, when \(D = D_{KL}\), we can write it as the
  difference between an (expected) energy and an entropy term:
  \begin{align*}
    \Fa(y)
    &= \E_{x \sim c'_\pi(y)} \left[ \log \frac{q(x|y)}{p_c(y|x) \cdot p_\pi(x)} \right] \\
    &= \E_{x \sim c'_\pi(y)} \left[ - \log p_c(y|x) - \log p_\pi (x) \right]
       - S_X \left[ c'_\pi(y) \right] \\
    &= \E_{x \sim c'_\pi(y)} \left[ E_{(\pi,c)}(x,y) \right] - S_X \left[ c'_\pi(y) \right]
  \; = U - TS
  \end{align*}
  where we call \(E_{(\pi,c)} : X \times Y \to \rr_+\) the \textit{energy} of
  the generative model \((\pi, c)\), and where \(S_X : \cat{C}(I, X) \to \rr_+\)
  is the Shannon entropy on \(X\).  The last equality makes the thermodynamic
  analogy: \(U\) is the \textit{internal energy} of the system; \(T = 1\) is the
  \textit{temperature}; and \(S\) is again the entropy.
\end{rmk}

Having now defined a more tractable fitness function, we can construct statistical games accordingly.
Since the free energy is an upper bound on relative entropy, optimizing the former can have the side effect of optimizing the latter\footnote{
Strictly speaking, one can have a decrease in free energy along with an increase in relative entropy, as long as the former remains greater than the latter.
Therefore, optimizing the free energy does not necessarily optimize the relative entropy.
However, as elaborated in Remark \ref{rmk:meaning-autoencoder}, the difference between the variational free energy and the relative entropy is the log-likelihood, so optimizing the free energy corresponds to simultaneous maximum-likelihood estimation and Bayesian inference.}.
We call the resulting games \textit{autoencoder games}, for reasons that will soon be clear.

\begin{ex}[Autoencoder] \label{ex:simp-autoenc-game}
  Let \(D : \cat{C}(I, X) \times \cat{C}(I, X) \to \rr\) be a measure of divergence between states on \(X\).
  Then a simple \(D\)\textit{-autoencoder game} is a simple statistical game \((c,\phi) : (X, X) \to (Y, Y)\) with fitness function \(\phi:\cctx(c)\to\rr\) given by
  \(\phi(\pi, k) = \E_{y\sim\efb{\pi}{c}{k}} \left[ \Fa_{D} \left(c'_\pi, c, \pi, y\right) \right]\)
  where \(c = (c, c') : (X,X) \lensto (Y,Y)\) constitutes the lens part of the game.
\end{ex}

One also of course has parameterized versions of the autoencoder games.

\begin{ex}[Simply parameterized autoencoder]
  A \textit{simply parameterized \(D\)-autoencoder game} is a simple parameterized statistical game \((c,\pi) : (X,X) \xto{(\Omega,\Omega)} (Y,Y)\) with the $D$-autoencoder fitness function \(\phi : \cctx(c) \to \rr\) given by
  \(\phi(\pi,k) = \E_{y\sim\efb{\pi}{c}{k}} \left[ \Fa_{D} \left(c'_\pi, c, \pi, y\right) \right]\). That is, a simply parameterized \(D\)-autoencoder game is just a simple \(D\)-autoencoder game with tensor product domain type.
\end{ex}

More often in applications, one doesn't use the same backwards channel to update both the ``belief about the causes'' in \(X\) and the parameters in \(\Omega\) simultaneously.
Instead, the backwards channel updates only the beliefs over \(X\), and any updating of the parameters is left to another `higher-order' process.
The \(X\)-update channel may nonetheless still itself be parameterized in \(\Omega\):
for instance, if it represents an approximate inference algorithm, then one often wants to be able to improve the approximation, and such improvement amounts to a change of parameters.
The `higher-order' process that performs the parameter updating is then often represented as a reparameterization:
a 2-cell in the bicategory \(\para(\SGame{\cat{C}})\).
The next example tells the first part of this story.

\begin{ex}[Parameterized autoencoder]
  A (simple) \textit{parameterized \(D\)-autoencoder game} is a parameterized statistical game \((c,\phi) : (X,X) \xto{(\Omega,\Theta)} (Y,Y)\) with fitness function \(\phi : \cctx(c) \to \rr\) given by
  \[
  \phi(\pi,k) = \E_{y\sim\efb{\pi}{c}{k}} \left[
    \Fa_{D} \bigl((c'_\pi)_X, c|{\pi_\Omega}, \pi_X, y\bigr)
  \right] \, .
  \]
  As before, the notation \((-)_X\) indicates taking the $X$ marginal, along the projection $\Theta\otimes X\klto X$.
  We also define \(c|\pi_\Omega := c\klcirc(\pi_\Omega \otimes \,\id_X)\), indicating ``$c$ given the parameter state $\pi_\Omega$''.
  Written out in full, the fitness function is therefore given by
  \[
  \phi(\pi,k) = \E_{y\sim\efb{\pi}{c}{k}} \left[
    \Fa_{D} \bigl(\mathsf{proj}_X \klcirc c'_\pi, c\klcirc(\pi_\Omega \otimes \,\id_X), \mathsf{proj}_X \klcirc \pi, y\bigr)
  \right] \, .
  \]
  Note that the pair \(\bigl(c|{\pi_\Omega}, (c'_{\pi_\Omega})_X\bigr)\) defines an unparameterized simple Bayesian lens \((X,X)\lensto(Y,Y)\), with \(c|{\pi_\Omega} : X\klto Y\) and \((c'_{\pi_\Omega})_X : Y\xklto{X} X\).
\end{ex}

\begin{rmk}[Meaning of `autoencoder'] \label{rmk:meaning-autoencoder}
  Why do we call autoencoder games thus? The name originates in machine
  learning, where one thinks of the forwards channel as `decoding' some latent
  state into a prediction of some generated data, and the backwards channel as
  `encoding' a latent state given an observation of the data; typically, the
  latent state space is thought to have lower dimensionality than the observed
  data space, justifying the use of this `compression' terminology. A slightly
  more precise way to see this is to consider an autoencoder game where the
  context and forwards channel are fixed. The only free variable available for
  optimization in the fitness function is then the backwards channel, and the
  optimum is obtained when the backwards channel equals the exact inversion of
  the forwards channel (given the prior in the context, and for all elements of
  the support of the state obtained from the continuation). Conversely, allowing
  only the forwards channel to vary, it is easy to see that the autoencoder
  fitness function is then equal to the fitness function of a maximum
  log-likelihood game (up to a constant). Consequently, optimizing the fitness
  of an autoencoder game in general corresponds to performing approximate
  Bayesian inference and maximum likelihood estimation simultaneously. The
  optimal strategy (lens or parameter) can then be considered as representing an
  `optimal' model of the process by which observations are generated, along with
  a recipe for inverting that model (and hence `encoding' the causes of the
  data). The prefix \textit{auto-} indicates that this model is learnt in an
  unsupervised manner, without requiring input about the `true' causes of the
  observations.
\end{rmk}

Some authors (in particular, \textcite{Knoblauch2019Generalized}) take a variant
of the \(D\)-autoencoder fitness function to define a generalization of Bayesian
inference: in an echo of Remark \ref{rmk:meaning-autoencoder}, the intuition
here is that Bayesian inference simply \textit{is} maximum likelihood
estimation, except `regularized' by the uncertainty encoded in the prior, which
stops the optimum strategy being trivially given by a Dirac delta
distribution. By allowing both the choice of likelihood function and divergence
measure to vary, one obtains a family of generalized inference
methods. Moreover, when one retains the standard choices of log-density as
likelihood and relative entropy as divergence, the resulting \textit{generalized
  Bayesian inference games} coincide with variational autoencoder games; then,
when the forwards channel (or its parameter) is fixed, both types of game
coincide with the Bayesian inference games of Example \ref{ex:bayes-game} above.

\begin{ex}[Generalized Bayesian inference \parencite{Knoblauch2019Generalized}] \label{ex:gen-bayes-game}
  Let \(D : \cat{C}(I, X) \times \cat{C}(I, X) \to \rr\) be a measure of divergence between states on \(X\), and let \(l : Y \otimes X \klto I\) be any effect on \(Y\otimes X\).
  Then a simple \textit{generalized \((l,D)\)-Bayesian inference game} is a simple statistical game \((c,\phi) : (X,X) \to (Y,Y)\) with fitness function \(\phi : \cctx(c) \to \rr\) given by
  \[
  \phi(\pi, k) = \E_{y \sim \efb{\pi}{c}{k}} \bigg[ \E_{x \sim c'_\pi(y)} \left[ l(y,x) \right] + D(c'_\pi(y), \pi) \bigg]
  \]
  where \((c, c') : (X,X) \lensto (Y,Y)\) constitutes the lens part of the game.
\end{ex}

\begin{prop}
  Generalized Bayesian inversion and autoencoder games coincide when \(D =
  D_{KL}\) and \(l = -\log p_c\), where \(p_c\) is a density function for the
  forwards channel \(c\).
  \begin{proof}
    Consider the \(D_{KL}\)-free energy. We have
    \begin{align*}
      \Fa_{D_{KL}} \left( c'_\pi, c, \pi, y \right)
      &= \E_{x\sim c'_\pi(y)} \left[ -\log p_c(y|x) - \log p_\pi(x) \right] - S_X\left[ c'_\pi(y) \right] \qquad\text{by Remark \ref{rmk:helmholtz}} \\
      &= \E_{x\sim c'_\pi(y)} \left[ -\log p_c(y|x) \right] + \E_{x\sim c'_\pi(y)} \left[ \log q(x|y) -\log p_\pi(x) \right] \\
      &= \E_{x\sim c'_\pi(y)} \left[ -\log p_c(y|x) \right] + D_{KL}\left(c'_\pi(y), \pi\right) \\
      &= \E_{x\sim c'_\pi(y)} \left[ l(y,x) \right] + D\left(c'_\pi(y), \pi\right)
    \end{align*}
    where \(q\) is a density function for \(c'_\pi\).
  \end{proof}
\end{prop}

Unsurprisingly, as in the autoencoder case, there are parameterized and simply
parameterized variants of generalized Bayesian inference games.

Finally, we remark that, in the case where $\cat{C} = \Cat{sfKrn}$, where $I$ is not terminal and morphisms into $I$ correspond to functions into $[0,\infty]$, composing a Bayesian lens and its context gives a lens $(I,I)\lensto(I,I)$:
both the forwards and backwards parts of this `$I$-endolens' return positive reals (which in this context Jacobs and colleagues call \textit{validities} \parencite{Jacobs2019Structured,Cho2015Introduction}), and which we can think of as ``the environment's measurements of its compatibility with the lens''.
In this case, we can therefore define \textit{validity games}, where the fitness function is simply given by computing the backwards validity\footnote{
Note that the backwards effect, as an $I$-state-dependent effect (or `vector'), already depends upon the forwards validity, so we do not need to include the forwards validity directly in the fitness computation.}.
Since such a fitness function measures the interaction of the lens with its environment, the corresponding statistical games may be of relevance in modelling multi-agent or otherwise interacting statistical systems---for instance, in modelling evolutionary dynamics.
We leave the exploration of this for future work.

\section{References}

\printbibliography[heading=none]

\appendix

\section{2-local contexts, graphically} \label{sec:apdx-2ctx}

To clarify the idea that the 2-local contexts for the factors of a tensor product game (or morphism more generally) are obtained by ``filling the hole'' on the left or right of the tensor, one can work in the monoidal bicategory of \(\Cat{V}\)-profunctors and use the associated graphical calculus \parencite{Roman2020Open} to depict the `hole' in the context and its filler.
In this section, we work with a general monoidal category \(\cat{C}\), which may or may not be a category of lenses or games.
Nonetheless, the basic idea is the same: a (complex) context for a morphism \(X \to Y\) is given by a triple of a residual denoted $\Theta$, a state \(I\to \Theta\otimes X\) and a `continuation' (or `effect') \(\Theta\otimes Y\to I\), coupled according to the coend quotient rule.

Below, we show how to obtain the object of right local 2-contexts for a tensor product morphism \(f\otimes f' : X\otimes X' \to Y\otimes Y'\), using the graphical calculus of \(\Cat{V\mdash Prof}\).
At each stage, we depict on the left the object named on the right. We start with a complex context for the tensor along with a `filler' object of morphisms \(X' \to Y'\), which is shown ``filling the (left-hand) hole''.
In the first step, we use the composition rule of \(\cat{C}\) to connect the matching `ports' on the domain \(X\).
We then couple the matching \(Y\) port using the coend and gather \(\Theta\) and \(Y\) together into a single residual.
Note that these steps correspond directly to factors in the definition of \(\pi_{f'} : \cctx(f\otimes f') \to \cctx(f')\) (Definition \ref{def:2local-ctx}).

\begin{gather*}
  \begin{tikzpicture}
	\begin{pgfonlayer}{nodelayer}
		\node [style=dot] (0) at (-7, 0) {$I$};
		\node [style=dot] (2) at (-3, 0) {$X$};
		\node [style=dot] (3) at (-3, -2) {$X'$};
		\node [style=dot] (6) at (3, 0) {$Y$};
		\node [style=dot] (7) at (-1, 0) {$X$};
		\node [style=dot] (8) at (1, 0) {$Y$};
		\node [style=dot] (9) at (3, -2) {$Y'$};
		\node [style=dot] (11) at (7, 0) {$I$};
		\node [style=action] (12) at (5, 0) {$\otimes$};
		\node [style=action] (13) at (-5, 0) {$\otimes$};
		\node [style=none] (14) at (-3, 2) {};
		\node [style=none] (15) at (3, 2) {};
	\end{pgfonlayer}
	\begin{pgfonlayer}{edgelayer}
		\draw (0) to (13);
		\draw [in=-180, out=0, looseness=1.25] (13) to (2);
		\draw [in=180, out=-90, looseness=1.25] (13) to (3);
		\draw [in=-180, out=0, looseness=1.25] (6) to (12);
		\draw (12) to (11);
		\draw [in=-90, out=0, looseness=1.25] (9) to (12);
		\draw (7) to (8);
		\draw [in=180, out=90, looseness=1.25] (13) to (14.center);
		\draw [in=360, out=90, looseness=1.25] (12) to (15.center);
		\draw (14.center) to (15.center);
	\end{pgfonlayer}
\end{tikzpicture} \qquad \int^{\Theta:\cat{C}} \cat{C}(I,\Theta\otimes X\otimes X')\times\cat{C}(X,Y)\times\cat{C}(\Theta\otimes Y\otimes Y',I) \\ \\
  \longrightarrow\quad \begin{tikzpicture}
	\begin{pgfonlayer}{nodelayer}
		\node [style=dot] (0) at (-5, 0) {$I$};
		\node [style=dot] (2) at (-1, -2) {$X'$};
		\node [style=dot] (3) at (1, 0) {$Y$};
		\node [style=dot] (5) at (-1, 0) {$Y$};
		\node [style=dot] (6) at (1, -2) {$Y'$};
		\node [style=dot] (7) at (5, 0) {$I$};
		\node [style=action] (8) at (3, 0) {$\otimes$};
		\node [style=action] (9) at (-3, 0) {$\otimes$};
		\node [style=none] (10) at (-1, 2) {};
		\node [style=none] (11) at (1, 2) {};
	\end{pgfonlayer}
	\begin{pgfonlayer}{edgelayer}
		\draw (0) to (9);
		\draw [in=180, out=-90, looseness=1.25] (9) to (2);
		\draw [in=-180, out=0, looseness=1.25] (3) to (8);
		\draw (8) to (7);
		\draw [in=-90, out=0, looseness=1.25] (6) to (8);
		\draw [in=180, out=90, looseness=1.25] (9) to (10.center);
		\draw [in=360, out=90, looseness=1.25] (8) to (11.center);
		\draw (10.center) to (11.center);
		\draw (9) to (5);
	\end{pgfonlayer}
\end{tikzpicture} \qquad\qquad\quad\quad \int^{\Theta:\cat{C}} \cat{C}(I,\Theta\otimes Y\otimes X')\times\cat{C}(\Theta\otimes Y\otimes Y',I) \\ \\
  \longhookrightarrow\quad \begin{tikzpicture}
	\begin{pgfonlayer}{nodelayer}
		\node [style=dot] (0) at (-5, 0) {$I$};
		\node [style=dot] (1) at (-1, -2) {$X'$};
		\node [style=dot] (4) at (1, -2) {$Y'$};
		\node [style=dot] (5) at (5, 0) {$I$};
		\node [style=action] (6) at (3, 0) {$\otimes$};
		\node [style=action] (7) at (-3, 0) {$\otimes$};
		\node [style=none] (8) at (-1, 2) {};
		\node [style=none] (9) at (1, 2) {};
	\end{pgfonlayer}
	\begin{pgfonlayer}{edgelayer}
		\draw (0) to (7);
		\draw [in=180, out=-90, looseness=1.25] (7) to (1);
		\draw (6) to (5);
		\draw [in=-90, out=0, looseness=1.25] (4) to (6);
		\draw [in=180, out=90, looseness=1.25] (7) to (8.center);
		\draw [in=360, out=90, looseness=1.25] (6) to (9.center);
		\draw (8.center) to (9.center);
	\end{pgfonlayer}
\end{tikzpicture} \qquad\qquad\qquad\qquad\; \int^{\Theta':\cat{C}} \cat{C}(I,\Theta'\otimes X')\times\cat{C}(\Theta'\otimes Y',I)
\end{gather*}

We find this graphical representation to be a useful aid in comprehension, and often simplifies the symbolic `book-keeping' that can complicate expressions such as those in Definition \ref{def:2local-ctx}.
The cost of this expressivity is the introduction of another categorical structure, and its associated cognitive load.
In previous work \parencite{Smithe2020Cyber}, we have made more use of this representation: there, we worked with the `optical' definition of Bayesian lenses described in \textcite{Smithe2020Bayesian};
and we note that an earlier informal version of this graphical language was originally used to define local contexts for tensor product games in the compositional game theory literature \parencite{Bolt2019Bayesian}.
Since the \(\Cat{V}\)-profunctorial setting \parencite{Clarke2020Profunctor} is required in order to define optics, we had already paid this extra cognitive cost.
In this paper, however, we have preferred to stick with the simpler fibrational definition of Bayesian lenses.

\begin{rmk}
  A different but related graphical calculus for optics is described by
  \textcite{Boisseau2020String}. However, this alternative calculus is somewhat
  less general than that of \textcite{Roman2020Open}, and its adoption here
  would not eliminate the extra cognitive cost; we do nonetheless make use of it
  in \parencite{Smithe2020Bayesian}.
\end{rmk}

\end{document}